\documentclass[11pt,a4paper]{amsart}
\usepackage[T1]{fontenc}
\usepackage[utf8]{inputenc}
\usepackage[textwidth=155mm,textheight=235mm,centering]{geometry}
\usepackage{amsmath,amssymb,amsthm,mathtools}
\usepackage{newtxtext,newtxmath}
\usepackage{microtype,booktabs,array,enumitem,needspace,float}
\usepackage{xcolor}
\definecolor{referenceblue}{RGB}{0,55,120}
\usepackage[colorlinks=true,linkcolor=referenceblue,citecolor=referenceblue,
  urlcolor=referenceblue,linktoc=all]{hyperref}
\usepackage{bookmark}
\DeclareRobustCommand{\namedref}[2]{\hyperref[#2]{#1~\ref*{#2}}}
\DeclareRobustCommand{\namedeqref}[2]{\hyperref[#2]{#1~\textup{(\ref*{#2})}}}
\renewcommand{\eqref}[1]{\hyperref[#1]{\textup{(\ref*{#1})}}}
\hypersetup{%
  pdftitle={m-Modified Conformal Vector Fields: An Adjoint Criterion and Four-Dimensional Classification},%
  pdfsubject={Pseudo-Riemannian geometry and metric Lie algebras},%
  pdfauthor={Li Zhigang; Deng Shaoqiang; Zhang Hui}%
}
\setlist[enumerate]{itemsep=3pt,topsep=5pt,leftmargin=2.4em}
\numberwithin{equation}{section}
\newtheorem{maintheorem}{Theorem}

\newtheorem{theorem}{Theorem}[section]
\newtheorem{proposition}[theorem]{Proposition}
\newtheorem{corollary}[theorem]{Corollary}
\newtheorem{lemma}[theorem]{Lemma}
\theoremstyle{definition}
\newtheorem{example}[theorem]{Example}
\theoremstyle{remark}
\newtheorem{remark}[theorem]{Remark}
\newtheorem*{remark*}{Remark}
\allowdisplaybreaks[2]
\newcommand{\R}{\mathbb R}
\newcommand{\C}{\mathbb C}
\newcommand{\Z}{\mathbb Z}
\newcommand{\g}{\mathfrak g}
\newcommand{\h}{\mathfrak h}
\newcommand{\nilradical}{\mathfrak n}
\newcommand{\ad}{\operatorname{ad}}
\newcommand{\Ad}{\operatorname{Ad}}
\newcommand{\tr}{\operatorname{tr}}
\newcommand{\im}{\operatorname{im}}
\newcommand{\rank}{\operatorname{rank}}
\newcommand{\spec}{\operatorname{spec}}
\newcommand{\Hess}{\operatorname{Hess}}
\newcommand{\grad}{\operatorname{grad}}
\newcommand{\Div}{\operatorname{div}}
\newcommand{\id}{\operatorname{Id}}
\newcommand{\aff}{\mathfrak{aff}}
\newcommand{\slTwoR}{\mathfrak{sl}_2(\R)}
\newcommand{\suTwo}{\mathfrak{su}(2)}
\newcommand{\heisenbergThree}{\mathfrak h_3}
\newcommand{\diag}{\operatorname{diag}}
\newcommand{\mMCVF}{\ensuremath{m}\textup{MCVF}}

\title[An adjoint criterion in even dimensions]{An Adjoint Criterion for Left-Invariant $m$-Modified Conformal Vector Fields on Even-Dimensional Lie Groups}
\author[Li]{Zhigang Li}
\author[Deng]{Shaoqiang Deng}
\author[Zhang]{Hui Zhang}
\date{}
\begin{document}
\begin{abstract}
  For a fixed positive real constant $m$, a real Lie algebra of even dimension $2k$ admits a
  nonzero left-invariant $m$-modified conformal vector field for some nondegenerate metric if
  and only if an adjoint operator has a nonzero real eigenvalue. In that case, a solution exists
  with neutral signature $(k,k)$ and zero conformal factor. As an application, we classify the
  admissible four-dimensional Lie algebras and metric signatures; the only Lorentzian exceptions
  among the admissible algebras are $\aff(\C)$ and the diamond algebra. On a pseudo-Riemannian
  manifold of arbitrary indefinite signature, a nonzero $m$-modified conformal vector field
  forces lightlike geodesic incompleteness, even when its conformal factor is variable.
\end{abstract}
\maketitle
\section{Introduction}

Which real Lie algebras admit a nonzero left-invariant $m$-modified conformal vector field for
some nondegenerate left-invariant metric? For fixed $m>0$, we characterize existence with
constant conformal factor in every even dimension and determine the admissible algebras and
metric signatures in dimension four.

On a pseudo-Riemannian manifold $(M,g)$ of dimension $n\ge2$, an \emph{$m$-modified conformal
vector field}, abbreviated to \mMCVF{}, is a vector field $X$ satisfying
\begin{equation}\label{eq:variable-conformal-equation}
  \frac12\mathcal L_Xg-\frac1m X^\flat\otimes X^\flat=\sigma g,
\end{equation}
where $X^\flat=g(X,\cdot)$ and $\sigma\in C^\infty(M)$ is the conformal factor. Unless stated
otherwise, we take $\sigma=c\in\R$ and consider
\begin{equation}\label{eq:constant-conformal-equation}
  \frac12\mathcal L_Xg-\frac1m X^\flat\otimes X^\flat=cg,
  \qquad c\in\R,
\end{equation}
including the case $c=0$.

Ordinary left-invariant conformal fields on unimodular metric Lie groups are Killing
\cite[Theorem 1.1]{AraujoChenLeandro2016}. Algebraic restrictions and classifications for
non-Killing left-invariant conformal fields in Lorentzian and trans-Lorentzian signatures are given in
\cite{ZhangChen2020Conformal,ZhangChen2022Lorentzian,ZhangChenTan2024}. Zhang and Chen classified
Riemannian Lie groups admitting nonzero left-invariant $m$-modified conformal vector fields
\cite[Theorem 4.2]{ZhangChen2023}. Their pseudo-Riemannian analysis assumes a nonzero conformal
factor \cite[Section 5]{ZhangChen2023}.
\namedeqref{Equation}{eq:constant-conformal-equation} is also a specialization of the
generalized Ricci soliton equation \cite{NurowskiRandall2016}. Nonzero solutions with $c=0$
already occur on groups locally isomorphic to $\mathrm{Sol}^3$ and $\mathrm{SL}_2(\R)$
\cite[Theorem 3.1(1)]{Calvaruso2017}.

Allowing the metric to vary, we obtain a criterion involving only the adjoint representation. We
write $\ad_Y(Z)=[Y,Z]$ and use $(p,q)$ for a signature with $p$ positive and $q$ negative
directions; $(k,k)$ is called neutral.

\Needspace{13\baselineskip}
\begin{maintheorem}\label{thm:adjoint-existence-criterion}
  Let $\g$ be a real Lie algebra of dimension $2k$, where $k\ge1$ is an integer, and fix $m>0$.
  The following conditions are equivalent:
  \begin{enumerate}[label=\textup{(\roman*)},ref=\roman*]
    \item\label{cond:nonzero-modified-conformal-field} Some nondegenerate symmetric form on $\g$
    admits a nonzero left-invariant \mMCVF{} with constant conformal factor.
    \item\label{cond:nonzero-real-adjoint-eigenvalue} $\ad_Y$ has a nonzero real eigenvalue for
    some $Y\in\g$.
    \item\label{cond:neutral-zero-factor-field} Some form of neutral signature $(k,k)$ on $\g$
    admits a nonzero left-invariant \mMCVF{} with zero conformal factor.
  \end{enumerate}
\end{maintheorem}

The proof constructs neutral metrics through explicit pairings on and between Jordan chains of
suitable nilpotent adjoint maps. The even-dimensional assumption is essential; see
\namedref{Example}{ex:odd-dimensional-counterexample}.

In dimension four, we apply the criterion to the classification and tabulations of real Lie
algebras in
\cite{AndradaEtAl2005,FreibertSchulteHengesbach2012,BiggsRemsing2016}. Write
$\R T\ltimes_D\mathfrak h$ for the extension by a derivation $D$, with $[T,v]=Dv$. Let
$\heisenbergThree=\operatorname{span}\{u,v,z\}$, with $[u,v]=z$. For $B\in M_2(\R)$, the
derivation $D_B$ acts as $B$ on $\operatorname{span}\{u,v\}$ and as $\tr B$ on $\R z$. The
algebras $\aff(\R)$ and $\aff(\C)$ are those of the affine groups of the real and complex lines.
A signature is Lorentzian if $\min\{p,q\}=1$.

\begin{maintheorem}\label{thm:four-dimensional-existence-and-signatures}
  Fix $m>0$. A four-dimensional real Lie algebra admits a nonzero left-invariant \mMCVF{} with
  constant conformal factor for some nondegenerate metric if and only if it is isomorphic to an
  algebra in the following list:
  \begin{enumerate}[label=\textup{\Roman*.},ref=\Roman*,leftmargin=*]
    \item\label{family:almost-abelian-overview} $\R T\ltimes_D\R^3$, with $D$ having a nonzero
    real eigenvalue;
    \item\label{family:heisenberg-extension-overview} $\R T\ltimes_{D_B}\heisenbergThree$, with
    $\tr B\ne0$ or $\det B<0$;
    \item\label{family:split-real-affine-overview} $\aff(\R)\oplus\aff(\R)$;
    \item\label{family:complex-affine-overview} $\aff(\C)$, regarded as a real Lie algebra;
    \item\label{family:split-simple-central-overview} $\slTwoR\oplus\R$.
  \end{enumerate}
  Every algebra in the list admits such a field for a metric of signature $(2,2)$. Among these
  algebras, a Lorentzian metric carrying such a field exists precisely for those other than
  $\aff(\C)$ and the diamond algebra
  \[
    \mathfrak d=\R T\ltimes_{D_{\diag(1,-1)}}\heisenbergThree.
  \]
  Every admissible indefinite signature can be realized with zero conformal factor. A Riemannian
  metric carrying such a field exists precisely for algebras in
  \namedref{family}{family:almost-abelian-overview} admitting a presentation in which $D$ is
  similar over $\R$ to
  \[
    aI_3+S,\qquad a\ne0,\qquad S^T=-S,
  \]
  where $I_3$ is the $3\times3$ identity matrix and $S$ is a real $3\times3$ matrix.
\end{maintheorem}

\namedref{Theorem}{thm:four-dimensional-existence-and-signatures} determines which
four-dimensional Lie algebras admit such fields and which metric signatures are possible. Its
Riemannian assertion follows from \cite[Theorem 4.2]{ZhangChen2023}. Different matrix parameters
may give isomorphic Lie algebras. We do not classify metrics up to automorphism.

Beyond the left-invariant setting, a nonzero \mMCVF{} in any indefinite signature forces
lightlike geodesic incompleteness, even for a variable conformal factor
(\namedref{Theorem}{thm:lightlike-incompleteness}). We also discuss compactness restrictions in
relation to the Riemannian triviality results in
\cite{PoddarSharma2025,DeshmukhEtAl2026ModifiedConformal}.

\namedref{Section}{sec:pseudo-riemannian-restrictions} establishes the geometric restrictions
and gives a compact quotient example. \namedref{Section}{sec:left-invariant-reduction} gives the
left-invariant algebraic reduction. \namedref{Section}{sec:adjoint-existence-criterion} proves
the adjoint criterion. \namedref{Section}{sec:four-dimensional-classification} gives the
four-dimensional Lie algebras, all zero-factor fields, and the admissible signatures.

\medskip
\noindent\textbf{Notation and conventions.}
All manifolds and tensor fields are smooth. Manifolds are connected, and vector fields are
globally defined unless a domain is specified. Throughout, $n\ge2$ denotes the dimension of the
manifold or Lie algebra, and $m>0$ is fixed. Metrics are nondegenerate and symmetric, with
arbitrary signature unless stated otherwise; a Riemannian metric is positive definite. The
signature $(p,q)$ counts positive and negative directions, respectively.

We use $c$ for the constant conformal factor and $\sigma$ for a factor that may vary, and write
$\alpha=X^\flat=g(X,\cdot)$ and $s=g(X,X)$. We write $\nabla$ for the Levi-Civita connection and
$\Hess f=\nabla df$. The metric adjoint $A^*$ is characterized by $g(AY,Z)=g(Y,A^*Z)$, whereas
$A^T$ denotes matrix transpose. We use $\id$ for the identity endomorphism and $I_r$ for the
$r\times r$ identity matrix, omitting the subscript when its size is clear.

A nonzero vector field means $X\not\equiv0$; a nowhere vanishing field satisfies $X_p\ne0$ at
every point. A vector $v$ is isotropic if $v\ne0$ and $g(v,v)=0$. For a vector field $V$,
isotropic means $g(V,V)=0$, with nonvanishing stated separately. A subspace $W$ is totally
isotropic if $g|_{W\times W}=0$. A lightlike geodesic is a nonconstant geodesic $\gamma$ with
$g(\dot\gamma,\dot\gamma)=0$. Lightlike geodesic completeness means that every maximal affinely
parametrized lightlike geodesic is defined on all of $\R$.

We use \emph{trans-Lorentzian} for signatures with $\min\{p,q\}=2$, equivalently $(n-2,2)$
up to overall sign, with $n\ge4$; compare \cite{ChenWolf2018,ZhangChenTan2024}. In dimension four,
neutral and trans-Lorentzian both mean signature $(2,2)$. By \namedref{Lemma}{lem:scaling}, replacing
$(g,X,c)$ by $(-g,-X,-c)$ preserves \eqref{eq:constant-conformal-equation}. Thus the existence
and vanishing statements for negative definite metrics follow from the Riemannian case.

An algebra is unimodular if $\tr\ad_Y=0$ for every $Y$. In explicit bases, unlisted brackets and
metric pairings are zero. Repeated indices do not imply summation; all sums over indices are
written explicitly.
\section{Pseudo-Riemannian geometric restrictions}\label{sec:pseudo-riemannian-restrictions}

We first record three differential identities following from the defining equation.

\begin{lemma}\label{lem:differential-identities}
  Let $X$ be an \mMCVF{} with conformal factor $\sigma\in C^\infty(M)$. Then
  \begin{equation}\label{eq:differential-identities}
    X(s)=2s\left(\sigma+\frac{s}{m}\right),\qquad
    \Div X=n\sigma+\frac{s}{m},\qquad
    \mathcal L_X\alpha=2\left(\sigma+\frac{s}{m}\right)\alpha.
  \end{equation}
  In the constant case, these identities apply with $\sigma=c$.
\end{lemma}
\begin{proof}
  Metric compatibility and the torsion-free identity $[X,Y]=\nabla_XY-\nabla_YX$ give
  \[
    \begin{aligned}
      (\mathcal L_Xg)(Y,Z)
      &=X\bigl(g(Y,Z)\bigr)-g([X,Y],Z)-g(Y,[X,Z])\\
      &=g(\nabla_YX,Z)+g(Y,\nabla_ZX).
    \end{aligned}
  \]
  Thus \eqref{eq:variable-conformal-equation} is equivalent to
  \begin{equation}\label{eq:symmetric-covariant-derivative}
    g(\nabla_YX,Z)+g(Y,\nabla_ZX)
    =2\sigma g(Y,Z)+\frac2m\alpha(Y)\alpha(Z).
  \end{equation}
  Setting $Y=Z=X$ and using metric compatibility once more yields
  \[
    X(s)=2g(\nabla_XX,X)=2\sigma s+\frac{2s^2}{m}.
  \]
  Write $\tr_g$ for contraction with the inverse metric. We have $\tr_g(\mathcal L_Xg)=2\Div X$,
  $\tr_g g=n$, and $\tr_g(\alpha\otimes\alpha)=g(X,X)=s$. Taking the metric trace of
  \eqref{eq:variable-conformal-equation} therefore gives
  \[
    \Div X-\frac{s}{m}=n\sigma,
  \]
  which is the second identity in \eqref{eq:differential-identities}. Finally, the derivation
  rule for the Lie derivative and $\mathcal L_XX=0$ give
  \[
    \mathcal L_X\alpha=(\mathcal L_Xg)(X,\cdot)
    =2\sigma\alpha+\frac2m\alpha(X)\alpha
    =2\left(\sigma+\frac{s}{m}\right)\alpha.
  \]
\end{proof}

The obstruction to completeness comes from the equation along a single lightlike geodesic. In
particular, it does not require a constant conformal factor. Lightlike geodesics were also used
to eliminate a multiple of the metric in \cite[proof of Theorem 1]{KiosakMatveev2009}. There,
in dimension $n>2$, if $g$ and $\widehat g=\psi^{-2}g$ are Einstein and $\psi>0$, then $\psi$
is affine along every affinely parametrized lightlike geodesic. For a global gradient solution
$X=\grad f$,
\eqref{eq:variable-conformal-equation} gives $\Hess f=\sigma g+m^{-1}df\otimes df$. Setting
$u=e^{-f/m}>0$ gives
\[
  \Hess u=-\frac{\sigma}{m}ug.
\]
Thus $u$ is affine along every affinely parametrized lightlike geodesic. The argument below
works directly with $g(X,\dot\gamma)$ and needs neither a potential nor an Einstein assumption.

\Needspace{8\baselineskip}
\begin{theorem}[Lightlike geodesic incompleteness]\label{thm:lightlike-incompleteness}
  Let $(M,g)$ be a pseudo-Riemannian manifold of signature $(p,q)$ with $p,q\ge1$. If it carries
  an \mMCVF{} $X\not\equiv0$, then it is lightlike geodesically incomplete. The conformal factor
  $\sigma\in C^\infty(M)$ need not be constant.
\end{theorem}
\begin{proof}
  Choose $o\in M$ with $X_o\ne0$. The isotropic vectors span $T_oM$. Indeed, choose a
  pseudo-orthonormal basis with positive vectors $e_i$ and negative vectors $v_j$. The vectors
  $e_i\pm v_1$ and $e_1\pm v_j$ are isotropic and span the whole space. Since $g$ is
  nondegenerate, there is an isotropic $w\in T_oM$ with $g(X_o,w)\ne0$. Replacing $w$ by $-w$ if
  necessary, arrange that $f_0=g(X_o,w)>0$. Let $\gamma$ be the maximal affinely parametrized
  geodesic with $\gamma(0)=o$ and $\dot\gamma(0)=w$. It is lightlike, since
  $g(\dot\gamma,\dot\gamma)$ is constant along a geodesic.

  Along this geodesic, set $f=g(X,\dot\gamma)$. Evaluating
  \eqref{eq:variable-conformal-equation} on $(\dot\gamma,\dot\gamma)$ and using
  $\nabla_{\dot\gamma}\dot\gamma=0$ and $g(\dot\gamma,\dot\gamma)=0$ gives
  \begin{equation}\label{eq:geodesic-riccati-equation}
    f'=g(\nabla_{\dot\gamma}X,\dot\gamma)=\frac{f^2}{m}.
  \end{equation}
  Since $f(0)=f_0>0$, the solution is
  \begin{equation}\label{eq:lightlike-riccati-solution}
    f(t)=\frac{f_0}{1-f_0t/m}
  \end{equation}
  throughout the geodesic's interval of definition containing $0$. If the geodesic were defined
  through $t=m/f_0$, the smooth function $g(X,\dot\gamma)$ would have the pole in
  \eqref{eq:lightlike-riccati-solution}. The maximal forward endpoint of $\gamma$ is therefore
  at most $m/f_0$, and in particular is finite.
\end{proof}

The following flat model realizes the pole in \eqref{eq:lightlike-riccati-solution} in every
indefinite signature and shows that the incompleteness mechanism does not require nonzero
curvature.

In the notation of \cite[Section 8, equations (88), (89), and (92)]{NurowskiRandall2016}, take
$c_1=1$, $c_2=\lambda=0$, $H=0$, and $f(x)=1/x$. With $x=r$, $y=u$, the solution one-form in
that reference is $-\alpha/m$ for the two-dimensional model below. Extra flat factors give the
stated higher-dimensional examples.

\Needspace{10\baselineskip}
\begin{example}\label{ex:flat-lightlike-incomplete-model}
  Let $n\ge2$. On the half-space $r<0$ with coordinates $(u,r,w_1,\ldots,w_{n-2})$, set
  \begin{equation}\label{eq:flat-left-invariant-model}
    g=2\,du\,dr+\sum_{j=1}^{n-2}\epsilon_j\,dw_j^2,
    \qquad X=-\frac mr\partial_u,\qquad \epsilon_j\in\{1,-1\},
  \end{equation}
  where $2\,du\,dr=du\otimes dr+dr\otimes du$. The metric is flat, and
  \[
    \alpha=-\frac mr\,dr,\qquad
    \frac12\mathcal L_Xg=\frac m{r^2}\,dr\otimes dr
    =\frac1m\alpha\otimes\alpha.
  \]
  Thus $c=0$. The substitution $r=-m e^{-t/m}$ gives
  \[
    X=e^{t/m}\partial_u,\qquad T=\partial_t,\qquad
    g=2e^{-t/m}\,du\,dt+\sum_{j=1}^{n-2}\epsilon_j\,dw_j^2.
  \]
  The frame $(T,X,\partial_{w_1},\ldots,\partial_{w_{n-2}})$ is left invariant for the product
  of $\R^{n-2}$ with the group
  \[
    (t,u)(s,v)=(t+s,u+e^{t/m}v).
  \]
  In this frame, $g(X,T)=1$, $g(T,T)=g(X,X)=0$, and the remaining pairings are constant, so $g$
  is also left invariant. Moreover, $[T,X]=X/m$. For each $n$, the choices of $\epsilon_j$
  realize every indefinite signature in dimension $n$. Fix $r_0<0$ and let
  \[
    \gamma(\tau)=(u_0,r_0+\tau,w_1^0,\ldots,w_{n-2}^0),
    \qquad \tau<-r_0.
  \]
  This is an affinely parametrized lightlike geodesic, and $g(X,\dot\gamma)=-m/(r_0+\tau)$ has a
  pole at $\tau=-r_0$. Thus the half-space is lightlike geodesically incomplete, and $X$ does
  not extend smoothly to the complete flat ambient space.
\end{example}

The Riemannian conclusion of the next theorem is given in \cite[Theorem 2.1]{PoddarSharma2025}.
The integrated identity in that proof also gives $c=s=0$ in arbitrary signature, using the
positive metric density. We give an alternative short proof using the extrema of the signed
squared length.

\begin{theorem}[Compact rigidity]\label{thm:compact-rigidity}
  Let $(M,g)$ be compact without boundary. Every \mMCVF{} with constant conformal factor $c$
  satisfies
  \[
    c=0,\qquad g(X,X)=0
    \quad\text{(hence $X=0$ if $g$ is Riemannian).}
  \]
  Its flow $\Phi_t$ is defined for all $t\in\R$ and obeys
  \begin{equation}\label{eq:compact-flow-metric}
    \Phi_t^*\alpha=\alpha,\qquad
    \Phi_t^*g=g+\frac{2t}{m}\alpha\otimes\alpha.
  \end{equation}
\end{theorem}
\begin{proof}
  Compactness ensures that $X$ has a complete flow $\Phi_t$. Apply
  \namedref{Lemma}{lem:differential-identities} with $\sigma=c$. At a maximum or minimum of
  $s=g(X,X)$, the first identity in \eqref{eq:differential-identities} gives $s\in\{0,-mc\}$. If
  $c>0$, this implies $-mc\le s\le0$. The second identity in \eqref{eq:differential-identities}
  then gives
  \[
    (n-1)c\le\Div X\le nc.
  \]
  If $c<0$, the same identity gives
  \[
    nc\le\Div X\le(n-1)c.
  \]
  Since $n\ge2$, in either case $\Div X$ has a strict constant sign, contradicting
  \[
    \int_M\Div X\,d\mu_g=0.
  \]
  Here $d\mu_g=|\det g|^{1/2}|dx^1\cdots dx^n|$ is the metric density in local coordinates. The
  identity $\mathcal L_Xd\mu_g=(\Div X)d\mu_g$ shows that the argument applies without
  orientability: differentiate $\int_M\Phi_t^*d\mu_g=\int_Md\mu_g$ at $t=0$. Thus $c=0$, and the
  same extremum argument gives $s\equiv0$.

  The third identity in \eqref{eq:differential-identities} now yields $\mathcal L_X\alpha=0$,
  while \eqref{eq:constant-conformal-equation} gives $\mathcal L_Xg=2\alpha\otimes\alpha/m$.
  Hence $\Phi_t^*\alpha=\alpha$ for all $t$, and
  \[
    \frac{d}{dt}\Phi_t^*g
    =\Phi_t^*(\mathcal L_Xg)
    =\frac{2}{m}\Phi_t^*(\alpha\otimes\alpha)
    =\frac{2}{m}\alpha\otimes\alpha,
  \]
  which gives \eqref{eq:compact-flow-metric} for every real $t$.
\end{proof}

For a solution in \namedref{Theorem}{thm:compact-rigidity}, we have $c=0$ and $g(X,X)=0$.
In Lorentzian signature, on the set where $X\ne0$, the identities
$\mathcal L_Xg=2\alpha\otimes\alpha/m$ and $\mathcal L_X\alpha=0$ are the constant-gauge
Kerr--Schild symmetry equations; see \cite[Section 2, Proposition 1]{CollEtAl2001}.
The flow formula \eqref{eq:compact-flow-metric} holds in every signature under the hypotheses
of that theorem.

When $\alpha$ is closed, the defining equation determines the full covariant derivative of $X$.
If $\alpha$ is globally exact, an exponential substitution converts the equation into a linear
Hessian equation.

\begin{proposition}\label{prop:closed-dual-and-gradient-reduction}
  An \mMCVF{} with constant conformal factor $c$ and $d\alpha=0$ satisfies
  \begin{equation}\label{eq:closed-dual-covariant-derivative}
    \nabla_YX=cY+\frac{\alpha(Y)}mX.
  \end{equation}
  If $\alpha=df$ globally, the function $u=e^{-f/m}>0$ satisfies
  \begin{equation}\label{eq:exponential-potential-hessian}
    \Hess u=-\frac cmu g.
  \end{equation}
  Conversely, every positive solution of \eqref{eq:exponential-potential-hessian} determines an
  \mMCVF{} by $X=-m\grad u/u$.

  If $X\not\equiv0$ and $s=g(X,X)$ is constant,
  \eqref{eq:closed-dual-covariant-derivative} also gives $c=-s/m$.
\end{proposition}
\begin{proof}
  Metric compatibility gives $(\nabla_Y\alpha)(Z)=g(\nabla_YX,Z)$. Together with the
  torsion-free identity, this implies
  \[
    d\alpha(Y,Z)=g(\nabla_YX,Z)-g(\nabla_ZX,Y),
  \]
  so $d\alpha=0$ makes $\nabla\alpha$ symmetric. Equation
  \eqref{eq:symmetric-covariant-derivative} with $\sigma=c$ therefore gives
  $\nabla\alpha=cg+m^{-1}\alpha\otimes\alpha$, or equivalently
  \eqref{eq:closed-dual-covariant-derivative}.

  If $\alpha=df$ and $u=e^{-f/m}$, then
  \[
    du=-\frac{u}{m}\alpha,
    \qquad
    \Hess u=\frac{u}{m^2}\alpha\otimes\alpha
    -\frac{u}{m}\nabla\alpha=-\frac cmu g.
  \]
  Conversely, if $u>0$ satisfies \eqref{eq:exponential-potential-hessian}, set $f=-m\log u$.
  Then
  \[
    \Hess f=-\frac m u\Hess u+\frac m{u^2}du\otimes du
    =cg+\frac1m df\otimes df,
  \]
  so $X=\grad f=-m\grad u/u$ satisfies \eqref{eq:constant-conformal-equation}.

  Finally, \eqref{eq:closed-dual-covariant-derivative} gives
  \[
    ds(Y)=2g(\nabla_YX,X)
    =2\left(c+\frac{s}{m}\right)\alpha(Y).
  \]
  If $s$ is constant and $X\not\equiv0$, evaluate this identity at a point where $X\ne0$.
  The covector $\alpha$ is nonzero there, so $c=-s/m$.
\end{proof}

The next conclusion is also a special case of \cite[Theorems 3--4]{MikesChodorova2010}. Those
results concern global potentials satisfying $\Hess f=\nu g+\tau\,df\otimes df$ on compact
manifolds without boundary, where $\nu$ and $\tau$ are smooth functions. The result for
indefinite metrics permits variable coefficients. In the Riemannian case, one additionally
requires $\nu\ge0$ everywhere or $\nu\le0$ everywhere. Here $\nu=c$ is constant, so this
condition holds. The constant conformal factor also lets us reduce directly to a parallel
gradient.

\begin{corollary}\label{cor:compact-global-gradient-vanishing}
  Every global gradient \mMCVF{} with constant conformal factor on a compact manifold without
  boundary vanishes identically, in any signature.
\end{corollary}
\begin{proof}
  Write $X=\grad f$ for a globally defined function $f$.
  \namedref{Theorem}{thm:compact-rigidity} gives $c=0$, so the positive function $u=e^{-f/m}$
  has parallel gradient by \eqref{eq:exponential-potential-hessian}. At a maximum of $u$ this
  gradient vanishes; parallel transport and connectedness imply that it vanishes everywhere.
  Hence $X=-m\grad u/u=0$.
\end{proof}

Closedness of $X^\flat$ alone is not sufficient for the conclusion of
\namedref{Corollary}{cor:compact-global-gradient-vanishing}. The following example gives compact
manifolds carrying a nonzero field whose dual form is closed but not exact.

The construction uses the product of the classical Lorentz--Sol geometry with a line. Compact
quotients and incompleteness of Lorentz--Sol geometry are known from
\cite{DumitrescuZeghib2010}. After a change of coordinates, the Lorentzian case below is an
instance of the product and compact quotient construction in
\cite[Examples 4.4--4.5]{HanounahMehidi2025}. The field with $c=0$ is obtained by
taking the $B=0$ case of \cite[Theorem 3.1(1)]{Calvaruso2017} and extending it by zero to the
extra factor. We compute the period of its dual form to distinguish closedness from global
exactness.

\Needspace{8\baselineskip}
\begin{example}\label{ex:compact-sol-quotient}
  There exist compact four-dimensional manifolds of signatures $(3,1)$ and $(2,2)$ carrying
  nonzero \mMCVF{}s $X$ with $c=0$. They are lightlike geodesically incomplete, and $X^\flat$ is
  closed but not exact.

  Consider the Lie group $G=(\R\ltimes_\Phi\R^2)\times\R$, with
  \[
    \Phi(t)=\diag(e^{t/m},e^{-t/m}),\qquad
    (t,w,r)(s,z,q)=(t+s,w+\Phi(t)z,r+q).
  \]
  Its underlying manifold is $\R\times\R^2\times\R$, so $G$ is simply connected. Write
  $w=(u,v)$, and set
  \[
    g=2e^{-t/m}\,dt\,du+e^{2t/m}\,dv^2+\epsilon\,dr^2,
    \qquad X=e^{t/m}\partial_u,\qquad \epsilon\in\{1,-1\}.
  \]
  Here $2\,dt\,du=dt\otimes du+du\otimes dt$. The group law shows that the coframe
  $(dt,e^{-t/m}du,e^{t/m}dv,dr)$ is left invariant.
  The coefficients of $g$ in this coframe are constant, and $X$ is the dual frame vector
  corresponding to $e^{-t/m}du$. Thus both $g$ and $X$ are left invariant. The $(t,u)$ block has
  negative determinant, giving signature $(3,1)$ when $\epsilon=1$ and $(2,2)$ when
  $\epsilon=-1$. Directly in these coordinates,
  \[
    \alpha=X^\flat=dt,\qquad
    \frac12\mathcal L_Xg=\frac1m\,dt\otimes dt
    =\frac1m\alpha\otimes\alpha.
  \]
  Thus $X$ satisfies \eqref{eq:constant-conformal-equation} with $c=0$.

  To obtain a lattice, put
  \[
    M=\begin{pmatrix}2&1\\1&1\end{pmatrix}\in\mathrm{SL}_2(\Z),\qquad
    \lambda=\frac{3+\sqrt5}{2},\qquad t_0=m\log\lambda.
  \]
  The eigenvalues of $M$ are $\lambda$ and $\lambda^{-1}$. Choose $P_0\in\mathrm{GL}_2(\R)$ with
  $\Phi(t_0)=P_0MP_0^{-1}$, and set $\Lambda=P_0\Z^2$. Since $\Phi(t_0)\Lambda=\Lambda$,
  \[
    \Gamma=\{(kt_0,\ell,j): k,j\in\Z,\ \ell\in\Lambda\}
  \]
  is a closed discrete subgroup. Its left action is free and properly discontinuous: for compact
  $K\subset G$, every $\gamma$ with $\gamma K\cap K\ne\varnothing$ lies in the compact set
  $KK^{-1}$, which meets $\Gamma$ in finitely many points. To see that $\Gamma\backslash G$ is
  compact, first reduce the $t$ coordinate modulo $t_0$, then the $\R^2$ coordinate modulo
  $\Lambda$, and finally the last coordinate modulo $\Z$. The product of the resulting closed
  fundamental regions is compact and maps onto the quotient. The left-invariant metric and field
  descend, and the descended field remains nonzero and satisfies the same equation. Since both
  signatures are indefinite, \namedref{Theorem}{thm:lightlike-incompleteness} gives lightlike
  geodesic incompleteness of the compact quotient.

  The descended form, again denoted by $\alpha$, is closed. The curve
  $\gamma(t)=\Gamma(t,0,0)$, $0\le t\le t_0$, is closed because $(t_0,0,0)\in\Gamma$, and
  \[
    \int_\gamma\alpha=\int_0^{t_0}dt=t_0\ne0.
  \]
  Thus $\alpha$ is not exact.
\end{example}

Here $\Gamma$ is not normal: a discrete normal subgroup of a connected Lie group is central, but
$(t_0,0,0)\in\Gamma$ acts nontrivially on $\R^2$ by conjugation. Thus $\Gamma\backslash G$ does
not inherit a quotient Lie group structure. This is consistent with the compact Lie-group
obstruction proved below in \namedref{Proposition}{prop:proper-flow-subgroup}.

The exponential substitution in \namedref{Proposition}{prop:closed-dual-and-gradient-reduction}
also has a geometric interpretation. Consider a nonzero left-invariant solution with $c=0$ and
$d\alpha=0$ on a simply connected Lie group with a left-invariant metric. Then $\alpha=df$
globally, and \eqref{eq:exponential-potential-hessian} shows that $\grad(e^{-f/m})=-e^{-f/m}X/m$
is parallel. This field is nowhere vanishing because $X$ is nonzero and left invariant. It is
isotropic because left invariance and \eqref{eq:differential-identities} give $g(X,X)=0$. By
\eqref{eq:closed-dual-covariant-derivative}, the field $X$ is recurrent: $\nabla_YX=\theta(Y)X$,
with $\theta=\alpha/m=d(f/m)$. The parallel field above is the standard rescaling of a recurrent
field with exact recurrence one-form; compare
\cite[pp.~1115--1116]{AlekseevskyCortesLeistner2023}.

We now turn to left-invariant fields, for which the defining equation becomes an algebraic
identity. This reduction will be used in the existence constructions that follow.

\section{Left-invariant fields in arbitrary signature}\label{sec:left-invariant-reduction}

Let $\g$ be the real Lie algebra of a connected Lie group. In the Lie-algebraic arguments below,
both $g$ and $X$ are left invariant. We identify $g$ with a nondegenerate symmetric form on
$\g$, and $X$ with a vector in $\g$. The equation becomes an identity of endomorphisms of $\g$.
We use the following notation:
\begin{equation}\label{eq:adjoint-rankone-decomposition}
  \begin{aligned}
    A&=\ad_X, &&\text{the adjoint action of }X,\\
    P&=\frac1mX\otimes\alpha, &&P(Y)=\frac{\alpha(Y)}mX,\\
    L&=A+c\id+P. &&
  \end{aligned}
\end{equation}
For left-invariant data, the conformal factor is necessarily constant. Indeed, evaluate the
equation on left-invariant $Y,Z$ with $g(Y,Z)\ne0$. All terms other than the conformal factor
are constant. The following proposition expresses the equation in terms of a skew-adjoint
endomorphism, in arbitrary signature.

\Needspace{14\baselineskip}
\begin{proposition}[Left-invariant reduction]\label{prop:left-invariant-skew-adjoint-reduction}
  For left-invariant $g$ and $X$, equation \eqref{eq:constant-conformal-equation} is equivalent
  to
  \begin{equation}\label{eq:metric-adjoint-field-equation}
    A+A^*=-2c\id-2P,
    \qquad\text{equivalently,}\qquad L^*=-L.
  \end{equation}
  Whenever these equivalent conditions hold, we have
  \begin{equation}\label{eq:trace-and-length-identities}
    s(s+mc)=0,\qquad
    \tr A=-nc-\frac{s}{m},\qquad
    A^*X=-2\left(c+\frac{s}{m}\right)X.
  \end{equation}
  In particular, $c\ne0$ implies that $\g$ is non-unimodular, whereas unimodularity implies
  $c=s=0$.
\end{proposition}
\begin{proof}
  For left-invariant $Y,Z$, the function $g(Y,Z)$ is constant, so the definition of the Lie
  derivative gives
  \[
    \begin{aligned}
      (\mathcal L_Xg)(Y,Z)
      &=X(g(Y,Z))-g([X,Y],Z)-g(Y,[X,Z])\\
      &=-g((A+A^*)Y,Z).
    \end{aligned}
  \]
  Also $g(PY,Z)=m^{-1}\alpha(Y)\alpha(Z)$, so $P^*=P$. Substitution into
  \eqref{eq:constant-conformal-equation} and nondegeneracy of $g$ give the first equation in
  \eqref{eq:metric-adjoint-field-equation}; \eqref{eq:adjoint-rankone-decomposition} then gives
  $L+L^*=A+A^*+2c\id+2P=0$. The first identity in \eqref{eq:trace-and-length-identities} follows
  from the first identity in \eqref{eq:differential-identities} with $\sigma=c$, because $s$ is
  constant. To obtain the other two, use $\tr P=\alpha(X)/m=s/m$, $\tr A^*=\tr A$, and $AX=0$:
  \[
    2\tr A=-2nc-\frac{2s}{m},\qquad
    A^*X=-2cX-2PX=-2\left(c+\frac{s}{m}\right)X.
  \]
  If $c\ne0$, the alternatives $s=0$ and $s=-mc$ give, respectively, $\tr A=-nc$ and
  $\tr A=-(n-1)c$, both nonzero.
\end{proof}

For later computations, write $G=(g(e_i,e_j))$ for the Gram matrix in a basis
$(e_1,\ldots,e_n)$, and let $\mathbf x$ be the coordinate column of $X$. We also write $A$ for
the matrix of $\ad_X$ in this basis. The coordinate column of $\alpha$ in the dual basis is
$G\mathbf x$. Thus \eqref{eq:metric-adjoint-field-equation} reads
\begin{equation}\label{eq:modified-conformal-gram-equation}
  A^TG+GA+2cG+\frac2m(G\mathbf x)(G\mathbf x)^T=0.
\end{equation}
With $Q=G^{-1}$, multiplication on the left and right by $Q$ gives the equivalent equation
\begin{equation}\label{eq:modified-conformal-inverse-gram-equation}
  AQ+QA^T+2cQ=-\frac2m\mathbf x\mathbf x^T.
\end{equation}
Thus, for fixed $X$ and $c$, the equation is linear in the inverse Gram matrix. Nondegeneracy
remains a separate requirement.

The following scaling identities also hold without the left-invariance assumption.

\begin{lemma}[Scaling]\label{lem:scaling}
  If $(g,X,c)$ satisfies \eqref{eq:constant-conformal-equation}, then so does $(kg,X/k,c/k)$ for
  every $k\in\R\setminus\{0\}$. If $m'>0$, then $(g,(m'/m)X,(m'/m)c)$ solves the corresponding
  equation with parameter $m'$. Thus the existence problem is independent of $m>0$, and
  replacing $(g,X,c)$ by $(-g,-X,-c)$ reverses the signature.
\end{lemma}
\begin{proof}
  In \eqref{eq:constant-conformal-equation}, the first transformation leaves both $\alpha$ and
  $\mathcal L_Xg$ unchanged, and also leaves $cg$ unchanged. For the second transformation, both
  sides are multiplied by $m'/m$, since $(m'/m)^2/m'=m'/m^2$.
\end{proof}

\Needspace{12\baselineskip}
The vanishing of the conformal factor imposes a stronger restriction than the trace identity
alone: the field $X$ lies in the image of its own adjoint map. This forces that adjoint map to
be nilpotent.

\begin{theorem}[Structure with zero conformal factor]\label{thm:zero-factor-structure}
  For a nonzero left-invariant \mMCVF{} with $c=0$,
  \begin{equation}\label{eq:zero-factor-identities}
    \begin{gathered}
      s=0,\qquad X\in\im A,\qquad A\ne0,\\
      A\text{ and }L\text{ are nilpotent},\qquad L^k=A^k\quad(k\ge2).
    \end{gathered}
  \end{equation}
  Conversely, among nonzero left-invariant \mMCVF{}s, the condition $X\in\im\ad_X$ implies
  $c=0$.
\end{theorem}
\begin{proof}
  By \eqref{eq:trace-and-length-identities}, $s=0$ and $A^*X=0$. It follows that
  \[
    P^2=AP=PA=0,
  \]
  where $AP=0$ uses $AX=0$, and $PA=0$ uses $\alpha\circ A=g(A^*X,\cdot)=0$.

  Suppose $X\notin\im A$. We have $(\im A)^\perp=\ker A^*$. Since $(W^\perp)^\perp=W$ for every
  subspace $W$, there is a $Y\in\ker A^*$ with $g(X,Y)\ne0$. Equation
  \eqref{eq:metric-adjoint-field-equation} gives
  \[
    AY=-2PY=-\frac{2}{m}g(X,Y)X.
  \]
  Since $g(X,Y)\ne0$, this gives $X\in\im A$, a contradiction. Hence $AZ=X$ for some $Z\in\g$,
  or equivalently $[Z,X]=-X$. The Jacobi identity gives $[\ad_Z,A]=\ad_{[Z,X]}=-A$.
  Consequently,
  \begin{equation}\label{eq:adjoint-dilation-identity}
    e^{t\ad_Z}Ae^{-t\ad_Z}=e^{-t}A\qquad(t\in\R),
  \end{equation}
  as follows by differentiating the left-hand side with respect to $t$. The spectrum of $A$ over
  $\C$ is consequently invariant under every positive dilation. Being finite, it consists only
  of zero, so $A$ is nilpotent. It is nonzero because $AZ=X\ne0$. The identities involving $P$
  now give $(A+P)^k=A^k$ for $k\ge2$, proving the assertions about $L=A+P$.

  For the converse, $X\in\im A$ gives \eqref{eq:adjoint-dilation-identity} by the same argument,
  independently of $c$, and hence $\tr A=0$. The alternatives in
  \eqref{eq:trace-and-length-identities} then give either $nc=0$ or $(n-1)c=0$. Since $n\ge2$,
  both imply $c=0$.
\end{proof}

\begin{corollary}\label{cor:nilpotent-lie-algebra-obstruction}
  A real nilpotent Lie algebra of dimension at least two admits no nonzero left-invariant
  \mMCVF{}, for any nondegenerate metric.
\end{corollary}
\begin{proof}
  Every adjoint map of a nilpotent Lie algebra is nilpotent. Thus $\g$ is unimodular, and
  \namedref{Proposition}{prop:left-invariant-skew-adjoint-reduction} forces $c=0$.
  \namedref{Theorem}{thm:zero-factor-structure} then supplies $Z$ with $\ad_ZX=-X$,
  contradicting nilpotence of $\ad_Z$.
\end{proof}

For comparison, Calvaruso \cite{Calvaruso2017} classifies left-invariant generalized Ricci
solitons on three-dimensional Riemannian and Lorentzian Lie groups. Jiang and Yang
\cite{JiangYang2025} study left-invariant generalized Ricci solitons on four-dimensional
nonabelian nilpotent Lie groups with left-invariant Lorentzian metrics.
\namedref{Corollary}{cor:nilpotent-lie-algebra-obstruction} gives a single
obstruction in every dimension $n\ge2$ and every signature, by the relation $[Z,X]=-X$, without
solving a metric-by-metric component system.

The following proposition gives a structural description and a uniform construction for
solutions with $c=0$ and $L=0$. The construction includes the Sol solution already present in
\cite[Theorem 3.1(1)]{Calvaruso2017}.

\Needspace{17\baselineskip}
\begin{proposition}\label{prop:rank-one-semidirect-construction}
  Let $\h$ be a finite-dimensional real Lie algebra containing a nonzero central vector $x$, and
  let $D$ be a derivation of $\h$ satisfying $Dx=\lambda x$ for some
  $\lambda\in\R\setminus\{0\}$. On $\g=\R T\ltimes_D\h$, where $[T,v]=Dv$, choose a vector-space
  decomposition $\h=\R x\oplus W$ and define
  \[
    g(x,x)=g(T,T)=0,\qquad g(x,T)=m\lambda,\qquad
    W\perp\operatorname{span}\{x,T\}.
  \]
  Give $W$ any nondegenerate symmetric form. Then $X=x$ is an \mMCVF{} with $c=0$ and $L=0$.

  Conversely, every nonzero left-invariant \mMCVF{} with $c=0$ and $L=0$ has a codimension-one
  ideal $\h=\ker\alpha$ containing $X$ in its center. There is an isotropic vector $T$ with
  $\alpha(T)=1$ such that
  \[
    \g=\R T\ltimes_D\h,\qquad DX=X/m.
  \]
\end{proposition}
\begin{proof}
  For the proposed metric, $\alpha(T)=m\lambda$ and $\alpha(\h)=0$. Since $x$ is central in
  $\h$, the map $A=\ad_x$ vanishes on $\h$ and satisfies $AT=-\lambda x$. Thus $A=-P$, and
  \namedref{Proposition}{prop:left-invariant-skew-adjoint-reduction} gives the required
  \mMCVF{}. The restriction of $g$ to $\operatorname{span}\{x,T\}$ has Gram matrix
  $\left(\begin{smallmatrix}0&m\lambda\\m\lambda&0\end{smallmatrix}\right)$, of determinant
  $-m^2\lambda^2<0$. This plane and $W$ are orthogonal and nondegenerate, so $g$ is
  nondegenerate.

  Conversely, \eqref{eq:adjoint-rankone-decomposition} gives $A=-P$, and
  \namedref{Theorem}{thm:zero-factor-structure} gives $s=0$. Applying the derivation identity
  for $A$ to $u,v\in\g$ gives
  \[
    -\frac{\alpha([u,v])}{m}X
    =-\frac{\alpha(u)}m[X,v]-\frac{\alpha(v)}m[u,X]=0.
  \]
  Indeed, $[X,v]=-\alpha(v)X/m$ and $[u,X]=\alpha(u)X/m$, so the two terms cancel. Since
  $X\ne0$, it follows that $\alpha([u,v])=0$. Thus $[\g,\g]\subset\ker\alpha$, so
  $\h=\ker\alpha$ is an ideal. Nondegeneracy of $g$ and $X\ne0$ give $\alpha\ne0$, hence $\h$
  has codimension one. Also $s=0$ gives $X\in\h$, and $A=-P$ gives $[X,\h]=0$. Choose $T_0$ with
  $\alpha(T_0)=1$ and put
  \[
    T=T_0-\frac12g(T_0,T_0)X.
  \]
  Then $\alpha(T)=1$ and $g(T,T)=0$. The derivation $D=\ad_T|_{\h}$ satisfies $DX=X/m$, giving
  the asserted decomposition.
\end{proof}

The algebraic restrictions also control the global behavior of the one-parameter subgroup.

\Needspace{8\baselineskip}
\begin{proposition}\label{prop:proper-flow-subgroup}
  Let $(G,g)$ be a connected Lie group of dimension at least two, with a left-invariant
  pseudo-Riemannian metric and a nonzero left-invariant \mMCVF{} $X$. Then $t\mapsto\exp(tX)$ is
  a proper embedding of $\R$, and its image is a closed subgroup of $G$. In particular, no
  compact Lie group admits such a field. No assumption of simple connectedness is needed.
\end{proposition}
\begin{proof}
  Write $\Ad:G\to\mathrm{GL}(\g)$ for the adjoint representation. We first show that
  \[
    \R\longrightarrow\mathrm{GL}(\g),\qquad
    t\longmapsto e^{tA}=\Ad(\exp(tX))
  \]
  is proper.

  If $c\ne0$, \namedref{Proposition}{prop:left-invariant-skew-adjoint-reduction} gives
  $\tr A\ne0$. On a compact subset of $\mathrm{GL}(\g)$, the absolute determinant is bounded
  above and away from zero. Since $\det(e^{tA})=e^{t\tr A}$, it follows that $t$ is bounded
  whenever $e^{tA}$ remains in such a subset.

  If $c=0$, \namedref{Theorem}{thm:zero-factor-structure} makes $A$ nonzero nilpotent. Let
  $r\ge1$ be maximal with $A^r\ne0$. Some entry of
  \[
    e^{tA}=\sum_{j=0}^{r}\frac{t^j}{j!}A^j
  \]
  is a polynomial of degree $r$ and tends to infinity in absolute value as $|t|\to\infty$. This
  proves properness in both cases.

  For compact $K\subset G$, the preimage of $K$ under $\exp(tX)$ is a closed subset of the
  compact preimage of $\Ad(K)$ under $e^{tA}$. Thus $t\mapsto\exp(tX)$ is proper as well. To
  prove injectivity, it suffices to show that $e^{tA}=I$ implies $t=0$. If $c\ne0$, this follows
  from the determinant. If $c=0$, factor
  \[
    e^{tA}-I=tA\sum_{j=0}^{r-1}\frac{(tA)^j}{(j+1)!}.
  \]
  The sum is the identity plus a nilpotent endomorphism and is therefore invertible. Since
  $A\ne0$, again $t=0$. The map $t\mapsto\exp(tX)$ is an immersion since its velocity is a
  translate of $X\ne0$. A proper injective immersion is a closed embedding, which also excludes
  compact $G$.
\end{proof}

\begin{remark*}[Compactness and real adjoint eigenvalues]
  For a compact Lie group $G$, the orbit $\Ad(G)X$ is compact for every $X\in\g$. If
  $\ad_YX=\lambda X$ for some $Y\in\g$, $X\ne0$, and $\lambda\in\R\setminus\{0\}$, then
  \[
    \Ad(\exp(tY))X=e^{t\ad_Y}X=e^{\lambda t}X\qquad(t\in\R).
  \]
  This orbit is unbounded as $t\to+\infty$ when $\lambda>0$ and as $t\to-\infty$ when
  $\lambda<0$, a contradiction. Thus, for every $Y\in\g$, the map $\ad_Y$ has no nonzero real
  eigenvalue.

  This gives another proof of the compact Lie-group obstruction in every dimension $n\ge2$.
  A nonzero left-invariant \mMCVF{} on a compact $G$ would have $c=0$ by
  \namedref{Theorem}{thm:compact-rigidity}, and
  \namedref{Theorem}{thm:zero-factor-structure} would supply $Z\in\g$ with $\ad_ZX=-X$,
  contradicting the preceding observation. In even dimensions, this also confirms the
  nonexistence predicted by \namedref{Theorem}{thm:adjoint-existence-criterion}: its
  real-eigenvalue condition fails for every $Y\in\g$. The compactness obstruction is valid for
  every $n\ge2$, whereas the equivalence in that theorem is asserted only in even dimensions.
  Noncompactness alone does not imply existence; for example, the noncompact abelian group
  $\R^{2k}$ has only zero adjoint maps.

  The compact manifold in \namedref{Example}{ex:compact-sol-quotient} illustrates the
  distinction between compact Lie groups and compact quotients by lattices: its metric and field
  descend from a noncompact group, and the quotient does not inherit a Lie group structure.
\end{remark*}

\section{The adjoint existence criterion}\label{sec:adjoint-existence-criterion}

We prove \namedref{Theorem}{thm:adjoint-existence-criterion} in this section.
\namedref{Proposition}{prop:even-dimensional-necessity} gives the necessary condition. For
sufficiency, assume that some adjoint map has a nonzero real eigenvalue. We consider two cases
to construct a neutral metric. First, suppose that a nonzero real eigenvalue occurs on the
nilradical, the largest nilpotent ideal of $\g$.
\namedref{Proposition}{prop:nilradical-real-eigenvalue-construction} then selects a vector that
satisfies the hypotheses of \namedref{Lemma}{lem:nilpotent-rank-one-metric-construction}. We
apply the lemma to construct an inverse Gram matrix. In the second case, no such eigenvalue
occurs on the nilradical. We then show that a noncompact simple ideal of a Levi subalgebra is a
direct summand of $\g$. Its Killing form makes the selected adjoint map skew-adjoint.
\namedref{Lemmas}{lem:nilpotent-skew-adjoint-signatures} and~\ref{lem:rank-one-metric-shear}
give the required neutral metric. The four-dimensional constructions and classification are
given in \namedref{Section}{sec:four-dimensional-classification}.

\Needspace{9\baselineskip}
\subsection{Necessary conditions and the role of even dimension}
\label{subsec:even-dimensional-necessity}

\begin{proposition}\label{prop:even-dimensional-necessity}
  In even dimension, a nonzero left-invariant \mMCVF{} can exist only if some adjoint map has a
  nonzero real eigenvalue.
\end{proposition}
\begin{proof}
  For $c=0$, \namedref{Theorem}{thm:zero-factor-structure} gives a vector $Z\in\g$ with
  $[X,Z]=X$, hence $\ad_ZX=-X$.

  If $c\ne0$ and $s=0$, then $A^*X=-2cX$ by \eqref{eq:trace-and-length-identities}. Since an
  endomorphism and its metric adjoint have the same characteristic polynomial, $-2c$ is also an
  eigenvalue of $A$.

  If $c\ne0$ and $s\ne0$, \eqref{eq:trace-and-length-identities} gives $s=-mc$ and $A^*X=0$.
  Hence $V=X^\perp$ is nondegenerate and $A$-invariant. Since $P|_V=0$,
  \eqref{eq:adjoint-rankone-decomposition} gives
  \[
    A|_V=-c\id+L|_V.
  \]
  By \eqref{eq:metric-adjoint-field-equation}, $L|_V$ is skew-adjoint on the odd-dimensional
  nondegenerate space $V$, so it is singular. Choose $0\ne v\in V$ with $Lv=0$. Then $Av=-cv$,
  which gives the required nonzero real eigenvalue.
\end{proof}

\begin{example}\label{ex:odd-dimensional-counterexample}
  The following three-dimensional example belongs to the Riemannian family in
  \cite[Theorem 4.2]{ZhangChen2023} and shows that the parity assumption cannot be omitted.
  Consider the Lie algebra
  \[
    [T,u]=au+bv,\qquad [T,v]=-bu+av,\qquad [u,v]=0,
    \qquad a,b\in\R,\quad ab\ne0.
  \]
  Equip it with the positive definite inner product for which $(T,u,v)$ is an orthonormal basis.
  Then $X=maT$ and $c=-ma^2$ satisfy \eqref{eq:constant-conformal-equation}. Indeed, on
  $\operatorname{span}\{u,v\}$ the symmetric part of $\ad_X$ is $ma^2\id$; on the $T$ direction
  the term $-\alpha\otimes\alpha/m$ equals $-ma^2g$, and all mixed components vanish. We write
  $\spec$ for the spectrum over $\C$, including algebraic multiplicities when listing
  eigenvalues. For $Y=tT+w$ with $w\in\operatorname{span}\{u,v\}$,
  \[
    \spec\ad_Y=\{0,t(a+ib),t(a-ib)\},
  \]
  so no adjoint map has a nonzero real eigenvalue.
\end{example}

\subsection{Neutral metric constructions and proof of the adjoint criterion}
\label{subsec:neutral-metric-constructions}

\Needspace{11\baselineskip}
\begin{lemma}[A neutral metric construction]\label{lem:nilpotent-rank-one-metric-construction}
  Let $A$ be a nilpotent endomorphism of a real vector space $V$ of dimension $2k$. Suppose that
  every even Jordan-block size at least four has even multiplicity, and let
  \[
    x\in(\im A\cap\ker A)\setminus\im A^2.
  \]
  Then, for every $m>0$, there is a nondegenerate symmetric form $g$ of signature $(k,k)$
  satisfying
  \begin{equation}\label{eq:rank-one-metric-identity}
    g(Au,v)+g(u,Av)=-\frac2m g(x,u)g(x,v)
    \qquad(u,v\in V).
  \end{equation}
\end{lemma}
\begin{proof}
  Fix $m>0$. Choose a Jordan basis whose first chain is $(e_1,e_2)$, with $e_1=x$ and
  $Ae_2=e_1$. Such a basis exists for the following reason. The vector $x$ is a linear
  combination of the bottom vectors of chains of length at least two. Since $x\notin\im A^2$, at
  least one coefficient from a chain of length two is nonzero. Use the coefficients in this
  linear combination on the second vectors of all the contributing chains to obtain $e_2$ with
  $Ae_2=x$. Replace that length-two chain by $(x,e_2)$ and leave all other chains unchanged. The
  coefficient on each of the two replaced basis vectors is the same nonzero number, so the
  replacement is invertible. Write $J_r$ for the Jordan block of size $r$, with ones on the
  superdiagonal, and let $\mathbf x=(1,0,\ldots,0)^T$.

  If $G$ is the Gram matrix of $g$ in this basis, then the required identity is
  \[
    A^TG+GA=-\frac2m(G\mathbf x)(G\mathbf x)^T.
  \]
  For $Q=G^{-1}$, this is equivalent to
  \begin{equation}\label{eq:rank-one-inverse-gram-equation}
    AQ+QA^T=-\frac2m\mathbf x\mathbf x^T.
  \end{equation}
  We construct a symmetric invertible solution $Q$ of neutral signature, block by block.

  For $r\ge1$, let $S_r$ be the signed antidiagonal matrix with $(S_r)_{i,r+1-i}=(-1)^{i-1}$.
  Direct calculation gives
  \[
    J_rS_r+S_rJ_r^T=0,\qquad
    S_r^T=(-1)^{r-1}S_r,\qquad S_rS_r^T=I_r.
  \]
  On a pair of equal Jordan blocks, the matrix
  \begin{equation}\label{eq:paired-jordan-inverse-gram-block}
    H_r=\begin{pmatrix}0&S_r\\S_r^T&0\end{pmatrix}
  \end{equation}
  satisfies the homogeneous equation
  \[
    \diag(J_r,J_r)H_r+H_r\diag(J_r^T,J_r^T)=0.
  \]
  It is nondegenerate and has an $r$-dimensional totally isotropic subspace, so its signature is
  $(r,r)$. For odd $r$, $S_r$ itself is symmetric and nondegenerate. Each antidiagonal pair
  contributes one positive and one negative direction; only the central entry is unpaired. Thus
  choosing $S_r$ or $-S_r$ gives either index difference $+1$ or $-1$, where the index
  difference means the positive index minus the negative index.

  Let $b\ge1$ be the number of Jordan blocks of size two. Set
  \[
    Q_x=-\frac1m\begin{pmatrix}0&1\\1&0\end{pmatrix}.
  \]
  This matrix has signature $(1,1)$ and satisfies
  \[
    J_2Q_x+Q_xJ_2^T=-\frac2m\diag(1,0).
  \]
  If $b$ is odd, assign $Q_x$ to the distinguished chain, and pair the remaining size-two blocks
  using $H_2$. If $b$ is even, pair the distinguished chain with one other size-two chain and
  use
  \[
    Q_{\mathrm{pair}}=\begin{pmatrix}Q_x&S_2\\S_2^T&0\end{pmatrix}.
  \]
  The homogeneous identity for $S_2$ gives
  \[
    \diag(J_2,J_2)Q_{\mathrm{pair}}
    +Q_{\mathrm{pair}}\diag(J_2^T,J_2^T)
    =-\frac2m\diag(1,0,0,0).
  \]
  The matrix $Q_{\mathrm{pair}}$ is nondegenerate because $S_2$ is invertible. Its last two
  coordinate directions span a totally isotropic subspace, so its signature is $(2,2)$. Pair
  the remaining $b-2$ size-two blocks using $H_2$.

  By hypothesis, the blocks of each even size $r\ge4$ can also be paired using $H_r$. Let $N$ be
  the total number of odd-sized Jordan blocks. Since $\dim V$ is even, $N$ is even. On half of
  these odd blocks choose the sign of $S_r$ to give index difference $+1$; on the other half
  choose it to give index difference $-1$.

  Take the direct sum of all these matrices, grouping paired chains together. Every summand is
  symmetric and nondegenerate. Only the distinguished summand contributes to the right-hand side
  of \eqref{eq:rank-one-inverse-gram-equation}. All even-dimensional summands are neutral, and
  the odd-dimensional summands have total index difference zero. Hence $Q$ has signature
  $(k,k)$. Its inverse $G=Q^{-1}$ has the same signature. Multiplication of
  \eqref{eq:rank-one-inverse-gram-equation} on the left and right by $G$ gives
  \eqref{eq:rank-one-metric-identity}, as required.

  More generally, the same construction realizes every signature $(k+j,k-j)$ with
  $j\in\mathbb Z$ and $|j|\le N/2$, by giving $N/2+j$ odd blocks index difference $+1$ and the
  others difference $-1$.
\end{proof}

\begin{proposition}[Sufficiency from a real adjoint eigenvalue on the nilradical]
  \label{prop:nilradical-real-eigenvalue-construction}
  Let $\g$ be a real Lie algebra of even dimension $2k$ and let $\nilradical$ be its nilradical.
  If $\ad_Y|_{\nilradical}$ has a nonzero real eigenvalue for some $Y\in\g$, then $\g$ admits a
  metric of signature $(k,k)$ and a nonzero left-invariant $m$-modified conformal vector field
  with zero conformal factor, for every $m>0$.
\end{proposition}
\begin{proof}
  Put $D=\ad_Y$, and write $\nilradical^1=\nilradical$ and
  $\nilradical^{j+1}=[\nilradical,\nilradical^j]$ for the lower central series. Each
  $\nilradical^j$ is an ideal of $\g$, hence is $D$-invariant. Choose $r$ maximal such that
  $D|_{\nilradical^r}$ has a nonzero real eigenvalue. Choose a corresponding real eigenvector
  $X\ne0$:
  \[
    X\in\nilradical^r,\qquad DX=\lambda X,\qquad
    \lambda\in\R\setminus\{0\}.
  \]
  By maximality, $D|_{\nilradical^{r+1}}$ has no nonzero real eigenvalue. In particular,
  $X\notin\nilradical^{r+1}$.

  Set $A=\ad_X$. Since $X\in\nilradical$, we have $A^j\g\subset\nilradical^j$ for every $j\ge1$,
  so $A$ is nilpotent. Moreover,
  \[
    AY=-\lambda X,\qquad AX=0,\qquad [D,A]=\lambda A.
  \]
  Thus $X\in\im A\cap\ker A$. Since $\nilradical^r$ is an ideal,
  \[
    A\g\subset\nilradical^r,\qquad
    A^2\g\subset[\nilradical^r,\nilradical^r]
    \subset\nilradical^{2r}\subset\nilradical^{r+1}.
  \]
  Consequently, $X\notin\im A^2$.

  It remains to check the required Jordan-block multiplicities. The identity
  $[D,A^j]=j\lambda A^j$ shows that every $\im A^j$ is $D$-invariant. Set
  \[
    E=\im A^2,\qquad F=\im A^3,
  \]
  and denote the restrictions of $D$ by $D_E$ and $D_F$. Since $E\subset\nilradical^{r+1}$, the
  only possible real eigenvalue of $D_E$ is zero. The inclusion $F\subset E$ gives
  $\spec D_F\subset\spec D_E$. On the other hand, $A:E\to F$ is surjective and satisfies
  \[
    D_FA=A(D_E+\lambda\id_E).
  \]
  The same identity shows that $\ker(A|_E)$ is $D_E$-invariant. The map $A$ induces an
  isomorphism $E/\ker(A|_E)\cong F$, under which the endomorphism induced by $D_E+\lambda\id_E$
  corresponds to $D_F$. It follows that
  \[
    \spec D_F\subset
    \spec D_E\cap\bigl(\spec D_E+\lambda\bigr).
  \]
  The two sets on the right can have real elements only at $0$ and $\lambda$, respectively.
  Since $\lambda\ne0$, $D_F$ has no real eigenvalue.

  For every $j\ge3$, the space $\im A^j$ is a $D$-invariant subspace of $F$, so $D$ has no real
  eigenvalue on it. Nonreal eigenvalues of a real endomorphism occur in conjugate pairs.
  Therefore $\rank A^j$ is even for every $j\ge3$. For a nilpotent map, the number of Jordan
  blocks of size $d\ge4$ is
  \[
    \rank A^{d-1}-2\rank A^d+\rank A^{d+1}.
  \]
  Indeed, a block of length $\ell$ contributes $\max\{\ell-j,0\}$ to $\rank A^j$. Hence every
  block size $d\ge4$ has even multiplicity.

  \namedref{Lemma}{lem:nilpotent-rank-one-metric-construction} now gives a metric of signature
  $(k,k)$ satisfying \eqref{eq:rank-one-metric-identity} with $x=X$. For left-invariant fields,
  \[
    (\mathcal L_Xg)(u,v)=-g(Au,v)-g(u,Av).
  \]
  Thus this metric and $X$ satisfy the modified conformal equation with zero conformal factor.
\end{proof}

\begin{lemma}[Signatures of forms making a nilpotent endomorphism skew-adjoint]
  \label{lem:nilpotent-skew-adjoint-signatures}
  Let $A$ be nilpotent on a real vector space of dimension $2k$. Suppose that $A$ is
  skew-adjoint for some nondegenerate symmetric form. Let $N$ be the number of odd-sized Jordan
  blocks of $A$, counted with multiplicity. Then $N$ is even. The signatures of nondegenerate
  symmetric forms making $A$ skew-adjoint are exactly
  \[
    (k+j,k-j),\qquad j\in\mathbb Z,\qquad |j|\le N/2.
  \]
  In particular, neutral signature $(k,k)$ is always possible.
\end{lemma}
\begin{proof}
  Let $V$ be the underlying space and let $B$ be a nondegenerate symmetric form for which $A$ is
  skew-adjoint. We construct a $B$-orthogonal decomposition into nondegenerate $A$-invariant
  subspaces. We use the fact that a totally isotropic subspace for a form of signature $(p,q)$
  has dimension at most $\min\{p,q\}$.

  At each step, work on a nonzero nondegenerate $A$-invariant subspace, using the same notation
  for the restrictions of $A$ and $B$. Let $r$ be the nilpotency index, so $A^r=0$ and
  $A^{r-1}\ne0$. The bilinear form
  \[
    C(u,v)=B(u,A^{r-1}v)
  \]
  is nonzero by nondegeneracy of $B$. For nonnegative integers $i,j$, skew-adjointness gives
  \[
    B(A^iu,A^jv)=(-1)^iB(u,A^{i+j}v),\qquad
    C(v,u)=(-1)^{r-1}C(u,v).
  \]
  Thus $C$ is symmetric when $r$ is odd and alternating when $r$ is even.

  Suppose first that $r=2q+1$ is odd. Then $C$ is symmetric. If $C(v,v)=0$ for every $v$, then,
  for all $u,w$,
  \[
    2C(u,w)=C(u+w,u+w)-C(u,u)-C(w,w)=0.
  \]
  This would imply $C=0$, a contradiction. Hence there is a vector $v$ with $\kappa=C(v,v)\ne0$.
  Put
  \[
    U=\operatorname{span}\{v,Av,\ldots,A^{r-1}v\}.
  \]
  For coefficients $a_0,\ldots,a_{r-1}$ not all zero, let $i$ be the smallest index with
  $a_i\ne0$ and set $w=\sum_{j=0}^{r-1}a_jA^jv$. The coefficients with $j<i$ vanish. The terms
  with $j>i$ involve a power of $A$ at least $r$, so
  \[
    \begin{aligned}
      B(w,A^{r-1-i}v)
      &=\sum_{j=i}^{r-1}(-1)^j a_j B(v,A^{j+r-1-i}v)\\
      &=(-1)^ia_i\kappa\ne0.
    \end{aligned}
  \]
  This proves both independence of the chain vectors and nondegeneracy of $B|_U$. Hence
  $\dim U=r$. The subspace $A^{q+1}U$ has dimension $q$ and is totally isotropic, since
  \[
    B(A^{q+1}x,A^{q+1}y)=(-1)^{q+1}B(x,A^{2q+2}y)=0
    \qquad(x,y\in U).
  \]
  Consequently $B|_U$ has signature $(q+1,q)$ or $(q,q+1)$.

  Now suppose that $r$ is even. Choose $u,v$ with $\kappa=C(u,v)\ne0$, and put
  \[
    W=\operatorname{span}\{A^ju,A^jv:0\le j\le r-1\}.
  \]
  For coefficients $(a_j,b_j)$ not all zero, let $i$ be the smallest index with
  $(a_i,b_i)\ne(0,0)$ and set $w=\sum_{j=0}^{r-1}(a_jA^ju+b_jA^jv)$. In the pairings below, the
  coefficients with $j<i$ are zero, and the terms with $j>i$ vanish because $A^r=0$. Since $C$
  is alternating, we obtain
  \[
    \begin{aligned}
      B(w,A^{r-1-i}v)
      &=\sum_{j=i}^{r-1}\bigl(a_jB(A^ju,A^{r-1-i}v)
      +b_jB(A^jv,A^{r-1-i}v)\bigr)\\
      &=(-1)^ia_i\kappa.
    \end{aligned}
  \]
  Similarly, using $C(v,u)=-\kappa$, we obtain
  \[
    B(w,A^{r-1-i}u)=(-1)^{i+1}b_i\kappa.
  \]
  At least one expression is nonzero. Thus the two chains are jointly independent and $B|_W$ is
  nondegenerate, with $\dim W=2r$. The subspace $A^{r/2}W$ has dimension $r$ and is totally
  isotropic, because
  \[
    B(A^{r/2}x,A^{r/2}y)=(-1)^{r/2}B(x,A^ry)=0
    \qquad(x,y\in W).
  \]
  It follows that $B|_W$ has signature $(r,r)$.

  Let $S$ be the subspace $U$ or $W$ just constructed. It is nondegenerate and $A$-invariant, so
  the current space is the orthogonal direct sum of $S$ and $S^\perp$. Moreover,
  \[
    B(Az,v)=-B(z,Av)=0\qquad(z\in S^\perp,\ v\in S),
  \]
  so $S^\perp$ is $A$-invariant. Its restricted form is nondegenerate, and the restricted
  operator remains nilpotent and skew-adjoint. Repeating the construction on this complement
  terminates by dimension. Each odd-dimensional summand contains one odd-sized Jordan block. Its
  positive and negative indices differ by one. Each even-dimensional summand contains two Jordan
  blocks of the same even size and has neutral signature. Thus there are exactly $N$
  odd-dimensional summands. Since $\dim V=2k$, their number $N$ is even.

  If $B$ has signature $(p,q)$, add the index differences of these summands:
  \[
    p+q=2k,\qquad p-q=\sum_{\ell=1}^{N}\varepsilon_\ell,
    \qquad\varepsilon_\ell\in\{1,-1\}.
  \]
  The sum is even and lies between $-N$ and $N$. Hence $(p,q)=(k+j,k-j)$ for an integer $j$ with
  $|j|\le N/2$. This proves that no other signatures are possible.

  Conversely, fix such an integer $j$. On each odd-dimensional summand, replacing the restricted
  form by its negative reverses its index difference and preserves skew-adjointness. Choose the
  signs so that $N/2+j$ summands have index difference $1$ and $N/2-j$ have index difference
  $-1$. Leave the forms on the even-dimensional summands unchanged. Their orthogonal direct sum
  is nondegenerate, has signature $(k+j,k-j)$, and still makes $A$ skew-adjoint. Taking $j=0$
  gives the neutral form.
\end{proof}

\begin{lemma}[A shear of a skew-adjoint metric]\label{lem:rank-one-metric-shear}
  Let $A$ be an endomorphism of a finite-dimensional real vector space, skew-adjoint for a
  nondegenerate symmetric form $B$. Suppose $AZ=X\ne0$ and $AX=0$. Set $\beta=B(X,\cdot)$ and
  $\zeta=B(Z,\cdot)$. Then
  \[
    g=B+\frac1m(\zeta\otimes\beta+\beta\otimes\zeta)
  \]
  is nondegenerate, has the same signature as $B$, and satisfies
  \eqref{eq:rank-one-metric-identity} with $x=X$.
\end{lemma}
\begin{proof}
  Skew-adjointness and the relations $AZ=X$, $AX=0$ give
  \[
    \begin{aligned}
      B(X,X)&=B(AZ,X)=-B(Z,AX)=0,\\
      B(X,Z)&=B(AZ,Z)=-B(Z,AZ)=-B(X,Z).
    \end{aligned}
  \]
  Thus $B(X,Z)=0$. Also, for every $u$,
  \[
    \beta(Au)=-B(AX,u)=0,\qquad
    \zeta(Au)=-B(AZ,u)=-\beta(u).
  \]
  The added tensor vanishes when one argument is $X$, so $g(X,\cdot)=\beta$. It follows that
  \[
    \begin{aligned}
      g(Au,v)+g(u,Av)
      &=\frac1m\bigl(\zeta(Au)\beta(v)+\beta(u)\zeta(Av)\bigr)\\
      &=-\frac2m\beta(u)\beta(v).
    \end{aligned}
  \]
  This proves \eqref{eq:rank-one-metric-identity}. Write $g(u,v)=B((\id+N)u,v)$, where
  \[
    N=\frac1m(Z\otimes\beta+X\otimes\zeta).
  \]
  This operator is self-adjoint for $B$. Since $NX=0$ and $NZ=B(Z,Z)X/m$, for every vector $u$
  we obtain
  \[
    N^2u=\frac{B(Z,Z)}{m^2}\beta(u)X,\qquad N^3u=0.
  \]
  In particular, no assumption on $B(Z,Z)$ is needed. Hence $(\id+tN)^{-1}=\id-tN+t^2N^2$ for
  every real $t$. The symmetric forms $B((\id+tN)\cdot,\cdot)$, $0\le t\le1$, provide a
  nondegenerate path from $B$ to $g$. The signature is therefore unchanged.
\end{proof}

\begin{proof}[Proof of \namedref{Theorem}{thm:adjoint-existence-criterion}]
  \namedref{Proposition}{prop:even-dimensional-necessity} proves
  \textup{(\ref{cond:nonzero-modified-conformal-field})}$\Rightarrow$\textup{(\ref{cond:nonzero-real-adjoint-eigenvalue})},
  and
  \textup{(\ref{cond:neutral-zero-factor-field})}$\Rightarrow$\textup{(\ref{cond:nonzero-modified-conformal-field})}
  is immediate. To prove
  \textup{(\ref{cond:nonzero-real-adjoint-eigenvalue})}$\Rightarrow$\textup{(\ref{cond:neutral-zero-factor-field})},
  let $\nilradical$ and $\mathfrak r$ be the nilradical and radical of $\g$. If some
  $\ad_Y|_{\nilradical}$ has a nonzero real eigenvalue, apply
  \namedref{Proposition}{prop:nilradical-real-eigenvalue-construction}.

  Suppose instead that, for every $Y\in\g$, the restriction $\ad_Y|_{\nilradical}$ has no
  nonzero real eigenvalue. Choose a Levi decomposition $\g=\mathfrak s\ltimes\mathfrak r$, with
  $\mathfrak s$ semisimple. The inclusion
  $[\g,\mathfrak r]\subset\nilradical$ gives
  \[
    \g/\nilradical\cong\mathfrak s\oplus(\mathfrak r/\nilradical),
  \]
  where the second summand is central and abelian.
  Relative to the invariant subspace $\nilradical$, each adjoint matrix is block triangular,
  with diagonal blocks given by the restriction and the quotient action. The nonzero real
  eigenvalue in \textup{(\ref{cond:nonzero-real-adjoint-eigenvalue})} must therefore occur on
  the quotient, hence on the adjoint representation of $\mathfrak s$. A compact semisimple
  algebra has a positive definite invariant form, so all its adjoint eigenvalues are purely
  imaginary. Thus $\mathfrak s$ has a noncompact simple factor $\mathfrak s_0$.

  We claim that $[\mathfrak s_0,\nilradical]=0$. A noncompact real simple algebra contains a
  subalgebra isomorphic to $\slTwoR$; see \cite[Proposition 6.52\textup{(b)}]{Knapp2002}. Choose
  a basis $(E,H,F)$ of this subalgebra satisfying
  \[
    [H,E]=2E,\qquad [H,F]=-2F,\qquad [E,F]=H.
  \]
  Consider the representation $\rho:\mathfrak s_0\to\mathfrak{gl}(\nilradical)$ defined by
  $\rho(U)v=[U,v]$. Suppose that $\rho$ is nonzero. Then $\ker\rho$ is a proper ideal of the
  simple algebra $\mathfrak s_0$, so $\ker\rho=0$. In particular, $\rho(H)\ne0$. Restrict $\rho$
  to the displayed $\slTwoR$ and complexify both the Lie algebra and its representation space,
  retaining the same notation. Finite-dimensional $\mathfrak{sl}_2(\C)$ representation theory
  implies that $\rho(H)$ is diagonalizable with integer eigenvalues. Since $\rho(H)\ne0$,
  at least one eigenvalue is nonzero. This is a nonzero
  real eigenvalue, contrary to our assumption. Thus $\rho=0$, which proves the claim.

  For $U,V\in\mathfrak s_0$ and $w\in\mathfrak r$, both $[V,w]$ and $[U,w]$ belong to
  $\nilradical$. Since $[\mathfrak s_0,\nilradical]=0$, Jacobi's identity gives
  \[
    [[U,V],w]=[U,[V,w]]-[V,[U,w]]=0.
  \]
  Every element of $\mathfrak s_0$ is a sum of brackets because
  $[\mathfrak s_0,\mathfrak s_0]=\mathfrak s_0$. Hence $[\mathfrak s_0,\mathfrak r]=0$. Let
  $\mathfrak s'$ be the sum of the other simple ideals of $\mathfrak s$ and put
  $\mathfrak h=\mathfrak s'+\mathfrak r$. The distinct simple factors commute, and $\mathfrak r$
  is an ideal, so $\mathfrak h$ is a subalgebra commuting with $\mathfrak s_0$. The Levi
  decomposition then gives
  \[
    \g=\mathfrak s_0\oplus\mathfrak h
  \]
  as a direct sum of ideals.

  Set $X=E$, $Z=-H/2$, and $A=\ad_X$. Then $AZ=X$ and $AX=0$. Since $[Z,X]=-X$,
  \namedeqref{equation}{eq:adjoint-dilation-identity} holds for this $Z$. The spectrum of $A$ is
  invariant under every positive dilation. Hence $A$ is nilpotent. By Cartan's criterion for
  semisimplicity, the Killing form on $\mathfrak s_0$ is
  nondegenerate. Its invariance makes $A|_{\mathfrak s_0}$ skew-adjoint. Extend it to
  $\g$ by taking an orthogonal direct sum with any nondegenerate symmetric form on
  $\mathfrak h$. Since $A|_{\mathfrak h}=0$, $A$ is skew-adjoint for this extended form.
  \namedref{Lemma}{lem:nilpotent-skew-adjoint-signatures} now supplies a form $B$ of signature
  $(k,k)$ for which $A$ is skew-adjoint. \namedref{Lemma}{lem:rank-one-metric-shear} then gives
  a metric of the same signature satisfying \eqref{eq:rank-one-metric-identity}. By
  \namedref{Proposition}{prop:left-invariant-skew-adjoint-reduction}, this is the required
  equation with zero conformal factor.
\end{proof}

\section{Four-dimensional classification}
\label{sec:four-dimensional-classification}

We first determine the Jordan types and signatures for a fixed zero-factor field. We then
classify the Lie algebras and all such fields. The final subsection determines which metric
signatures occur without assuming $c=0$.

\subsection{Jordan types and signatures}
\label{subsec:four-dimensional-jordan-signatures}

In dimension four, the possible Jordan chains determine the admissible signatures for a fixed
field with zero conformal factor. Chains of length two use the inverse Gram matrices from
\namedref{Lemma}{lem:nilpotent-rank-one-metric-construction}; a chain of length three uses
\namedref{Lemma}{lem:rank-one-metric-shear}. We first determine the possible Jordan types of
$\ad_X$.

\begin{lemma}\label{lem:four-dimensional-jordan-chains}
  Suppose that $\dim\g=4$ and that $0\ne X\in\im\ad_X$. Then $A=\ad_X$ is nilpotent. There is a
  basis $(e_1,e_2,e_3,e_4)$ with $e_1=X$ in which the nonzero images of the basis vectors are
  given by one of the following rows:
  \begin{equation}\label{eq:four-dimensional-jordan-chains}
    \begin{array}{c|l}
      \text{Jordan partition}&\text{nonzero }Ae_i\\ \hline
      (2,1,1)&Ae_2=e_1,\\
      (3,1)&Ae_2=e_1,\quad Ae_3=e_2,\\
      (2,2)&Ae_2=e_1,\quad Ae_4=e_3.
    \end{array}
  \end{equation}
\end{lemma}
\begin{proof}
  Choose $Z$ with $AZ=X$. The dilation identity \eqref{eq:adjoint-dilation-identity} and its
  spectral argument show that $A$ is nilpotent. This argument uses only $[Z,X]=-X$. Since
  $AZ=X\ne0$, the possible Jordan partitions are $(2,1,1)$, $(3,1)$, $(2,2)$, and $(4)$.

  Suppose the Jordan partition is $(4)$. In this case, $\ker A=\R X$, so choose a Jordan chain
  with $e_1=X$ and $Ae_i=e_{i-1}$ for $i=2,3,4$. Repeated use of the derivation rule gives
  \[
    \begin{aligned}
      0=A^4[e_3,e_4]
      &=\sum_{j=0}^4\binom4j[A^je_3,A^{4-j}e_4]\\
      &=4[e_2,e_1]+6[e_1,e_2]
      =2e_1.
    \end{aligned}
  \]
  Here $[e_1,e_2]=Ae_2=e_1$, which contradicts $e_1=X\ne0$.

  For $(2,1,1)$ and $(3,1)$, the space $\ker A\cap\im A$ is the bottom line of the unique
  nontrivial chain. Rescale that chain so that its bottom vector is $X$.

  For $(2,2)$, choose a basis $(X,e_3)$ of $\ker A=\im A$ and choose preimages $Ae_2=X$,
  $Ae_4=e_3$. These choices give the stated bases.
\end{proof}

\Needspace{13\baselineskip}
\begin{proposition}\label{prop:four-dimensional-jordan-signatures}
  Suppose that $\dim\g=4$ and $0\ne X\in\im\ad_X$, and fix $m>0$. The signatures of
  nondegenerate symmetric forms for which this fixed $X$ is an \mMCVF{} with zero conformal
  factor are exactly as follows:
  \begin{center}
    \normalfont
    \begin{tabular}{@{}ll@{}}
      \toprule
      Jordan partition of $\ad_X$ & Admissible signatures\\
      \midrule
      $(2,1,1)$ & $(3,1)$, $(2,2)$, $(1,3)$\\
      $(3,1)$ & $(3,1)$, $(2,2)$, $(1,3)$\\
      $(2,2)$ & $(2,2)$\\
      \bottomrule
    \end{tabular}
  \end{center}
\end{proposition}
\begin{proof}
  By \namedref{Proposition}{prop:left-invariant-skew-adjoint-reduction}, the field equation with
  $c=0$ is \eqref{eq:rank-one-metric-identity}, with $A=\ad_X$ and $x=X$. The first identity in
  \eqref{eq:trace-and-length-identities} gives $g(X,X)=0$. Since $X\ne0$, neither $(4,0)$ nor
  $(0,4)$ is possible. Thus every admissible signature is indefinite. Choose an adapted basis
  from \namedref{Lemma}{lem:four-dimensional-jordan-chains}, and write $(e^1,e^2,e^3,e^4)$ for
  its dual basis. The matrices $Q_x$ and $S_r$ are those in the proof of
  \namedref{Lemma}{lem:nilpotent-rank-one-metric-construction}. For each Jordan type, we also
  record a neutral Gram matrix for use in the explicit constructions below. The subscripts on
  these Gram matrices indicate the Jordan partition of $A$.

  Suppose the Jordan partition is $(2,1,1)$. Assign $Q_x$ to the distinguished chain. The two
  one-dimensional blocks may have independent signs. The inverse Gram matrices
  \[
    Q=\diag(Q_x,\varepsilon_1,\varepsilon_2),\qquad
    Q_x=-\frac1m\begin{pmatrix}0&1\\1&0\end{pmatrix},
    \qquad \varepsilon_1,\varepsilon_2\in\{1,-1\},
  \]
  satisfy \eqref{eq:rank-one-inverse-gram-equation} by the construction in
  \namedref{Lemma}{lem:nilpotent-rank-one-metric-construction}. The block $Q_x$ has signature
  $(1,1)$, so the sign choices give all three indefinite signatures. Inverting preserves the
  signature. For $(\varepsilon_1,\varepsilon_2)=(1,-1)$, a neutral representative is
  \begin{equation}\label{eq:neutral-gram-single-length-two-chain}
    G_{(2,1,1)}=
    \begin{pmatrix}
      0&-m&0&0\\
      -m&0&0&0\\
      0&0&1&0\\
      0&0&0&-1
    \end{pmatrix}.
  \end{equation}

  Suppose the Jordan partition is $(3,1)$. Here $X=A^2e_3$, so the hypothesis $X\notin\im A^2$
  in \namedref{Lemma}{lem:nilpotent-rank-one-metric-construction} fails. Instead, let $B_0$ have
  Gram matrix $\diag(-S_3,-1)$. Since $S_3^{-1}=S_3$, the homogeneous block identity in that
  proof shows that $A$ is skew-adjoint for $B_0$. There are two odd-sized Jordan blocks. Thus
  \namedref{Lemma}{lem:nilpotent-skew-adjoint-signatures}, with $k=2$ and $N=2$, gives forms of
  signatures $(3,1)$, $(2,2)$, and $(1,3)$ for which $A$ is skew-adjoint. Apply
  \namedref{Lemma}{lem:rank-one-metric-shear} to each form with $Z=e_2+e_1/(2m)$. Since $AZ=X$
  and $AX=0$, this gives solutions with the same respective signatures and the same $X$. For the
  neutral representative $B_0$ above, we have
  \[
    \beta=B_0(X,\cdot)=-e^3,\qquad
    \zeta=B_0(Z,\cdot)=e^2-\frac1{2m}e^3.
  \]
  The resulting metric is
  \[
    g_{(3,1)}=B_0
    -\frac1m(e^2\otimes e^3+e^3\otimes e^2)
    +\frac1{m^2}e^3\otimes e^3.
  \]
  Its Gram matrix is
  \begin{equation}\label{eq:neutral-gram-length-three-chain}
    G_{(3,1)}=
    \begin{pmatrix}
      0&0&-1&0\\
      0&1&-1/m&0\\
      -1&-1/m&1/m^2&0\\
      0&0&0&-1
    \end{pmatrix}.
  \end{equation}
  The $m^{-2}$ entry comes from the chosen preimage $Z$ of $X$.

  \Needspace{13\baselineskip}
  Suppose the Jordan partition is $(2,2)$. Pair the distinguished chain with the other chain,
  using the pairing block $-S_2$ in the construction of
  \namedref{Lemma}{lem:nilpotent-rank-one-metric-construction}. This sign change preserves the
  construction. It gives
  \[
    Q_{(2,2)}=\begin{pmatrix}Q_x&-S_2\\-S_2^T&0\end{pmatrix},\qquad
    S_2=\begin{pmatrix}0&1\\-1&0\end{pmatrix},
  \]
  whose inverse is
  \begin{equation}\label{eq:neutral-gram-paired-length-two-chains}
    G_{(2,2)}=
    \begin{pmatrix}
      0&0&0&-1\\
      0&0&1&0\\
      0&1&0&-1/m\\
      -1&0&-1/m&0
    \end{pmatrix}.
  \end{equation}
  The cited construction gives signature $(2,2)$ and \eqref{eq:rank-one-metric-identity}. To
  exclude other signatures in this case, suppose $g$ satisfies
  \eqref{eq:rank-one-metric-identity} and put $K=\ker A=\im A$. For $u\in K$, applying that
  identity first to $(u,u)$ and then to $(u,v)$ gives
  \[
    0=-\frac2m g(X,u)^2,\qquad
    g(u,Av)=-\frac2m g(X,u)g(X,v)=0
    \quad(v\in\g).
  \]
  Thus $K$ is a two-dimensional totally isotropic subspace. If $g$ has signature $(p,q)$, then
  \[
    2=\dim K\le\min\{p,q\},\qquad p+q=4.
  \]
  Hence $p=q=2$, which excludes both Lorentzian signatures.
\end{proof}

\begin{remark}\label{rem:four-dimensional-metric-construction}
  In dimension four, the construction starts with any eigenvector $X$ of $\ad_Y$ for a nonzero
  real eigenvalue and then uses the Jordan form of $\ad_X$. It therefore determines existence on
  the Lie algebra without requiring a normal form for its full bracket. The matrices
  \eqref{eq:neutral-gram-single-length-two-chain}--\eqref{eq:neutral-gram-paired-length-two-chains}
  are representatives establishing existence; these matrices do not classify all metrics
  admitting such fields up to automorphism.
\end{remark}

\subsection{Lie algebras and zero-factor fields}
\label{subsec:four-dimensional-lie-algebras}

We use the solvable classification in \cite[Theorem 1.5]{AndradaEtAl2005} and the following
nonsolvable classification from \cite[Section 2.2]{BiggsRemsing2016}.

\begin{lemma}\label{lem:four-dimensional-nonsolvable-classification}
  A four-dimensional nonsolvable real Lie algebra is isomorphic to $\suTwo\oplus\R$ or
  $\slTwoR\oplus\R$. In particular, no four-dimensional real Lie algebra is semisimple.
\end{lemma}

We now apply \namedref{Theorem}{thm:adjoint-existence-criterion} to these classifications, using
the notation for semidirect products, $\heisenbergThree$, and $D_B$ fixed in the introduction.

\begin{theorem}\label{thm:four-dimensional-lie-algebra-classification}
  A four-dimensional real Lie algebra admits a nonzero left-invariant \mMCVF{} for some
  nondegenerate metric if and only if it is isomorphic to one of the following:
  \begin{enumerate}[label=\textup{\Roman*.},ref=\Roman*,leftmargin=*]
    \item\label{family:almost-abelian} $\R T\ltimes_D\R^3$, where $\R^3$ is abelian and
    $D\in M_3(\R)$ has a nonzero real eigenvalue;
    \item\label{family:heisenberg-extension} $\R T\ltimes_{D_B}\heisenbergThree$, where
    \begin{equation}\label{eq:heisenberg-extension-eigenvalue-condition}
      \tr B\ne0\quad\text{or}\quad\det B<0;
    \end{equation}
    \item\label{family:split-real-affine} $\aff(\R)\oplus\aff(\R)$, with $[t_1,u_1]=u_1$ and
    $[t_2,u_2]=u_2$;
    \item\label{family:complex-affine} $\aff(\C)$, regarded as a real Lie algebra, with basis
    $(t,j,u,v)$ and
    \begin{equation}\label{eq:complex-affine-lie-brackets}
      [t,u]=u,\quad [t,v]=v,\quad [j,u]=v,\quad [j,v]=-u;
    \end{equation}
    \item\label{family:split-simple-central} $\slTwoR\oplus\R$.
  \end{enumerate}
  Every algebra in the list admits a trans-Lorentzian metric carrying a nonzero left-invariant
  \mMCVF{} with $c=0$. Different matrix parameters in \namedref{families}{family:almost-abelian}
  and~\ref{family:heisenberg-extension} may yield isomorphic Lie algebras.
\end{theorem}

\begin{proof}
  By \namedref{Theorem}{thm:adjoint-existence-criterion}, we must determine which algebras have
  an adjoint map with a nonzero real eigenvalue. The classification in \cite[Theorem
  1.5]{AndradaEtAl2005} gives every four-dimensional solvable real Lie algebra in one of the
  forms \namedref{families}{family:almost-abelian}--\ref{family:complex-affine}, before imposing
  the spectral restrictions in \namedref{families}{family:almost-abelian}
  and~\ref{family:heisenberg-extension}.

  Suppose $\g=\R T\ltimes_D\R^3$. For $Y=tT+w$, with $t\in\R$ and $w\in\R^3$, the matrix of
  $\ad_Y$ relative to $\R T\oplus\R^3$ is block triangular, with diagonal blocks $0$ and $tD$.
  Thus some adjoint map has a nonzero real eigenvalue if and only if $D$ does.

  Suppose $\g=\R T\ltimes_{D_B}\heisenbergThree$. For $Y=tT+w$, with $t\in\R$ and
  $w\in\heisenbergThree$, the invariant flag $\R z\subset\heisenbergThree\subset\g$ gives
  \begin{equation}\label{eq:heisenberg-extension-characteristic-polynomial}
    \det(\xi\id-\ad_{tT+w})
    =\xi(\xi-t\tr B)\det(\xi I_2-tB).
  \end{equation}
  If $\tr B\ne0$, then $\ad_T$ has the nonzero real eigenvalue $\tr B$. If $\tr B=0$, the
  eigenvalues of $B$ satisfy $\lambda^2=-\det B$. Hence
  \eqref{eq:heisenberg-extension-eigenvalue-condition} is precisely the required condition.

  In \namedref{families}{family:split-real-affine} and~\ref{family:complex-affine}, the maps
  $\ad_{t_1}$ and $\ad_t$, respectively, have eigenvalue $1$.

  Suppose $\g$ is nonsolvable. By
  \namedref{Lemma}{lem:four-dimensional-nonsolvable-classification}, the possibilities are
  $\suTwo\oplus\R$ and $\slTwoR\oplus\R$. On $\suTwo\oplus\R$, all adjoint maps are skew-adjoint
  for a common positive definite inner product, so none has a nonzero real eigenvalue. On
  $\slTwoR\oplus\R$, the standard triple satisfies $[H,E]=2E$.

  \namedref{Theorem}{thm:adjoint-existence-criterion} now gives the list and a metric of
  signature $(2,2)$ with $c=0$ for every algebra in it.
\end{proof}

\Needspace{10\baselineskip}
\begin{corollary}\label{cor:four-dimensional-excluded-algebras}
  Up to isomorphism, the four-dimensional real Lie algebras admitting no nonzero left-invariant
  \mMCVF{} are:
  \begin{enumerate}[label=\textup{(\roman*)},ref=\roman*]
    \item the nilpotent algebras;
    \item $\suTwo\oplus\R$;
    \item the oscillator algebra $[T,u]=v$, $[T,v]=-u$, $[u,v]=z$;
    \item\label{case:excluded-rotation-extension} $\R\oplus\mathfrak r'_{3,a}$, $a\ge0$, where
    $[T,u]=au+v$, $[T,v]=-u+av$ on $\mathfrak r'_{3,a}$.
  \end{enumerate}
\end{corollary}
\begin{proof}
  Apply \namedref{Theorem}{thm:four-dimensional-lie-algebra-classification} to the normal forms
  in \cite[Theorem 1.5]{AndradaEtAl2005}.

  Suppose $\g$ is solvable but not nilpotent. The excluded normal forms are
  $\R\oplus\mathfrak r'_{3,a}$, $a\ge0$, and the oscillator algebra, denoted by
  $\mathfrak d'_{4,0}$ in that classification. Their parameters satisfy, respectively,
  \[
    \spec D=\{0,a+i,a-i\},\qquad
    \tr B=0,\quad\det B=1.
  \]
  Neither satisfies the corresponding condition in
  \namedref{Theorem}{thm:four-dimensional-lie-algebra-classification}. Every other nonnilpotent
  solvable normal form satisfies one of its existence conditions.

  Nilpotent algebras are excluded by
  \namedref{Corollary}{cor:nilpotent-lie-algebra-obstruction}.

  Among the nonsolvable algebras,
  \namedref{Theorem}{thm:four-dimensional-lie-algebra-classification} excludes precisely
  $\suTwo\oplus\R$.
\end{proof}

\begin{remark}
  For example, $\R\oplus\mathfrak r'_{3,1}$ is non-unimodular but admits no nonzero
  left-invariant \mMCVF{}. Non-unimodularity is therefore not sufficient for existence. By
  \namedref{Proposition}{prop:left-invariant-skew-adjoint-reduction}, existence for some
  left-invariant metric depends only on the Lie algebra. A metric and a vector at the identity
  extend by left translations to any connected Lie group with that Lie algebra; simple
  connectivity is not required.
\end{remark}

\Needspace{12\baselineskip}
\begin{proposition}[All zero-factor fields in dimension four]
  \label{prop:four-dimensional-zero-factor-fields}
  Fix $m>0$ and a presentation in
  \namedref{Theorem}{thm:four-dimensional-lie-algebra-classification}.
  \namedref{Table}{tab:four-dimensional-zero-factor-fields} lists all nonzero vectors $X$ that
  are left-invariant \mMCVF{}s with $c=0$ for some nondegenerate symmetric metric, together with
  the Jordan partition of $\ad_X$ and all admissible signatures for that $X$. Every listed
  signature is realizable for every vector in its row. The metric is allowed to depend on $X$.
\end{proposition}

Use the bases in \namedref{Theorem}{thm:four-dimensional-lie-algebra-classification}. In
\namedref{family}{family:heisenberg-extension}, put $V=\operatorname{span}\{u,v\}$. In
\namedref{family}{family:split-simple-central}, choose a standard triple $(H,E,F)$ with
$[H,E]=2E$, $[H,F]=-2F$, and $[E,F]=H$. All coefficients in the table are real, and
$\lambda\in\R\setminus\{0\}$.

\begin{table}[H]
  \centering
  \small
  \renewcommand{\arraystretch}{1.22}
  \begin{tabular}{@{}p{14mm}p{65mm}p{21mm}p{38mm}@{}}
    \toprule
    Family & All admissible nonzero $X$ & \shortstack[l]{Jordan\\partition} & Admissible
    signatures\\
    \midrule
    \ref{family:almost-abelian}
    & $X=x\in\R^3\setminus\{0\}$, $Dx=\lambda x$
    & $(2,1,1)$ & $(3,1),(2,2),(1,3)$\\
    \addlinespace
    \ref{family:heisenberg-extension}
    & $X=rz$, $r\ne0$, with $\tr B\ne0$
    & $(2,1,1)$ & $(3,1),(2,2),(1,3)$\\
    \ref{family:heisenberg-extension}
    & $X=w+rz$, $w\in V\setminus\{0\}$, $Bw=\lambda w$, $r\in\R$
    & $(2,2)$ & $(2,2)$\\
    \addlinespace
    \ref{family:split-real-affine}
    & $X=pu_1+qu_2$, $(p,q)\ne(0,0)$, $pq=0$
    & $(2,1,1)$ & $(3,1),(2,2),(1,3)$\\
    \ref{family:split-real-affine}
    & $X=pu_1+qu_2$, $pq\ne0$
    & $(2,2)$ & $(2,2)$\\
    \addlinespace
    \ref{family:complex-affine}
    & $X=pu+qv$, $(p,q)\ne(0,0)$
    & $(2,2)$ & $(2,2)$\\
    \addlinespace
    \ref{family:split-simple-central}
    & $X=aH+bE+dF\ne0$, $a^2+bd=0$
    & $(3,1)$ & $(3,1),(2,2),(1,3)$\\
    \bottomrule
  \end{tabular}
  \caption{All nonzero left-invariant zero-factor fields when the metric
  may vary. The family parameters satisfy
  \namedref{Theorem}{thm:four-dimensional-lie-algebra-classification}.
  For each family represented by two rows, the sets of vectors in those
  rows are disjoint.}
  \label{tab:four-dimensional-zero-factor-fields}
\end{table}

\begin{proof}
  For a fixed nonzero $X$, \namedref{Theorem}{thm:zero-factor-structure} and
  \namedref{Proposition}{prop:four-dimensional-jordan-signatures} give the following criterion:
  a metric with $c=0$ exists if and only if $X\in\im\ad_X$. Set $A=\ad_X$. Every admissible
  vector therefore satisfies
  \[
    X\in\im A\subset[\g,\g].
  \]
  We determine all such vectors and the Jordan partitions of their adjoint maps.

  Suppose $\g$ belongs to \namedref{family}{family:almost-abelian}. Since $[\g,\g]\subset\R^3$,
  write $X=x\in\R^3\setminus\{0\}$. The brackets give
  \[
    AT=-Dx,\qquad A|_{\R^3}=0,\qquad \im A=\R Dx.
  \]
  Thus $0\ne X\in\im A$ is equivalent to $Dx=\lambda x$ for some nonzero real $\lambda$.
  Conversely, this equation gives $A(-T/\lambda)=X$. Here $\rank A=1$ and $A^2=0$, so the
  partition is $(2,1,1)$.

  Suppose $\g$ belongs to \namedref{family}{family:heisenberg-extension}. Since
  $[\g,\g]\subset V+\R z$, write $X=w+rz$, with $w\in V$. Define the alternating form $\omega$
  on $V$ by $[y_1,y_2]=\omega(y_1,y_2)z$ for $y_1,y_2\in V$; thus $\omega(u,v)=1$. The brackets
  give
  \[
    AT=-Bw-r(\tr B)z,\qquad Ay=\omega(w,y)z\quad(y\in V),
    \qquad Az=0.
  \]

  Suppose $w=0$. Since $X\ne0$, we have $r\ne0$. The condition $X\in\im A$ holds exactly when
  $\tr B\ne0$. In this case $\rank A=1$ and $A^2=0$, giving $(2,1,1)$.

  Suppose $w\ne0$. Nondegeneracy of $\omega$ gives $A(V)=\R z$, and hence
  \[
    \im A=\operatorname{span}\{Bw,z\}.
  \]
  Thus $X=w+rz$ lies in $\im A$ exactly when $Bw=\lambda w$ for a nonzero real $\lambda$; the
  coefficient $r$ is unrestricted. For every such vector,
  \[
    \rank A=2,\qquad A^2T=-\omega(w,Bw)z=0,\qquad
    A^2|_{V+\R z}=0.
  \]
  The Jordan partition is therefore $(2,2)$.

  Suppose $\g$ belongs to \namedref{family}{family:split-real-affine}. The inclusion
  $X\in\im A\subset[\g,\g]$ gives $X=pu_1+qu_2$, with $(p,q)\ne(0,0)$. For every such vector,
  \[
    At_1=-pu_1,\qquad At_2=-qu_2,\qquad Au_1=Au_2=0,
    \qquad A(-t_1-t_2)=X.
  \]
  Thus all these vectors satisfy the criterion. If $pq=0$, the rank is one and the partition is
  $(2,1,1)$. If $pq\ne0$, the rank is two and the partition is $(2,2)$. In both cases $A^2=0$.

  Suppose $\g$ belongs to \namedref{family}{family:complex-affine}. Again $X\in[\g,\g]$, so
  $X=pu+qv$ with $(p,q)\ne(0,0)$. The brackets give
  \[
    At=-X,\qquad Aj=qu-pv,\qquad Au=Av=0.
  \]
  The two displayed images are independent because $p^2+q^2\ne0$. Every such $X$ therefore lies
  in $\im A$, and $\rank A=2$, $A^2=0$. This gives $(2,2)$.

  Suppose $\g$ belongs to \namedref{family}{family:split-simple-central}. Since
  $X\in[\g,\g]=\slTwoR$, its additional central component is zero. Write $X=aH+bE+dF$. In the
  standard two-dimensional representation,
  \[
    X=\begin{pmatrix}a&b\\d&-a\end{pmatrix},\qquad
    X^2=(a^2+bd)I_2.
  \]
  If $X\in\im\ad_X$, then $X=[X,Z]$ for some $Z\in\slTwoR$. Cyclicity of the matrix trace gives
  \[
    2(a^2+bd)=\tr(X^2)=\tr\bigl(X[X,Z]\bigr)
    =\tr(X^2Z)-\tr(XZX)=0.
  \]
  Conversely, suppose $a^2+bd=0$, so $X^2=0$. For $X\ne0$, choose $v_0\in\R^2$ with $Xv_0\ne0$.
  The vectors $Xv_0,v_0$ are independent, and in this basis $X$ has matrix
  $E=\left(\begin{smallmatrix}0&1\\0&0\end{smallmatrix}\right)$. This change of basis induces an
  automorphism of $\slTwoR$. Now
  \[
    [E,-H/2]=E,\qquad \ad_EH=-2E,\qquad \ad_EF=H,\qquad \ad_EE=0.
  \]
  It follows that $X\in\im\ad_X$ and that its adjoint map has one block of size three, together
  with the central block of size one. The partition is $(3,1)$.

  We have proved both necessity and sufficiency for every row. The signature column follows from
  \namedref{Proposition}{prop:four-dimensional-jordan-signatures}.
\end{proof}

For example, on $\aff(\R)\oplus\aff(\R)$, the vector $u_1$ can be realized with all three
indefinite signatures, whereas $u_1+u_2$ can be realized only with signature $(2,2)$. Thus the
signatures available for a particular vector can be more restricted than those available
somewhere on the same Lie algebra.

\subsection{Signature classification}\label{subsec:four-dimensional-signatures}

\namedref{Proposition}{prop:four-dimensional-zero-factor-fields} classifies all zero-factor
vectors and their signatures. We now determine which signatures occur on each Lie algebra
without assuming $c=0$. The notation $\sim_{\R}$ denotes similarity over $\R$. Let $\mathfrak d$
denote the real diamond algebra with basis $(T,u,v,z)$ and nonzero brackets
\begin{equation}\label{eq:diamond-lie-brackets}
  [T,u]=u,\qquad [T,v]=-v,\qquad [u,v]=z.
\end{equation}
Every algebra in \namedref{family}{family:heisenberg-extension} with $\tr B=0$ and $\det B<0$ is
isomorphic to $\mathfrak d$. Indeed, $B$ has eigenvalues $\lambda,-\lambda$ with $\lambda>0$.
Choose corresponding eigenvectors $u,v$, replace $T$ by $T/\lambda$, and set $z=[u,v]$.

\Needspace{16\baselineskip}
\begin{theorem}\label{thm:four-dimensional-signature-classification}
  Let $\g$ be a four-dimensional real Lie algebra and let $m>0$. The following statements
  concern nonzero left-invariant \mMCVF{}s.
  \begin{enumerate}[label=\textup{(\roman*)},ref=\roman*]
    \item\label{signature:neutral-existence} An \mMCVF{} for a trans-Lorentzian metric exists
    precisely for the algebras in
    \namedref{Theorem}{thm:four-dimensional-lie-algebra-classification}.
    \item\label{signature:lorentzian-existence} An \mMCVF{} for a Lorentzian metric exists
    precisely for the algebras in that list other than $\aff(\C)$ and $\mathfrak d$.
    \item\label{signature:riemannian-existence} An \mMCVF{} for a Riemannian metric exists
    precisely for the algebras with a presentation
    \begin{equation}\label{eq:riemannian-extension-derivation}
      \begin{gathered}
        \g=\R T\ltimes_D\R^3,\qquad D\sim_{\R}aI_3+S,\\
        a\in\R\setminus\{0\},\qquad S\in M_3(\R),\qquad S^T=-S.
      \end{gathered}
    \end{equation}
    The condition on $D$ is necessary and sufficient for each given presentation in
    \namedref{family}{family:almost-abelian}.
  \end{enumerate}
  All the indefinite signatures admitted in \textup{(\ref{signature:neutral-existence})} and
  \textup{(\ref{signature:lorentzian-existence})} can be realized with $c=0$.
\end{theorem}

The zero-factor vectors and their signatures are listed in
\namedref{Table}{tab:four-dimensional-zero-factor-fields}. For $c=0$,
\eqref{eq:trace-and-length-identities} forces $g(X,X)=0$, so no Riemannian metric admits a
nonzero solution with $c=0$.

Signature reversal is understood as in \namedref{Lemma}{lem:scaling}; in particular, the
Riemannian results also determine existence for negative definite metrics.

\Needspace{8\baselineskip}
\begin{proof}[Proof of
\namedref{Theorem}{thm:four-dimensional-signature-classification}\textup{(\ref{signature:neutral-existence})--(\ref{signature:lorentzian-existence})}]
  Part~\textup{(\ref{signature:neutral-existence})} follows from
  \namedref{Theorem}{thm:four-dimensional-lie-algebra-classification}. By
  \namedref{Lemma}{lem:scaling}, it suffices to consider signature $(3,1)$ in
  part~\textup{(\ref{signature:lorentzian-existence})}.
  \namedref{Proposition}{prop:four-dimensional-zero-factor-fields} and
  \namedref{Table}{tab:four-dimensional-zero-factor-fields} give a metric of this signature with
  $c=0$ for every algebra in \namedref{families}{family:almost-abelian},
  \ref{family:split-real-affine}, and~\ref{family:split-simple-central}, and for every algebra
  in \namedref{family}{family:heisenberg-extension} with $\tr B\ne0$.

  It remains to exclude Lorentzian solutions on $\mathfrak d$ and $\aff(\C)$.

  Suppose $\g=\mathfrak d$. This algebra is unimodular, so
  \namedref{Proposition}{prop:left-invariant-skew-adjoint-reduction} gives $c=0$. By
  \namedref{Proposition}{prop:four-dimensional-zero-factor-fields}, $\ad_X$ has Jordan partition
  $(2,2)$ for every admissible $X$, so only neutral signature is possible. Thus no Lorentzian
  solution exists.

  Suppose $\g=\aff(\C)$. The same proposition excludes Lorentzian solutions with $c=0$. It
  remains to exclude $c\ne0$. Use the basis $(t,j,u,v)$ of
  \eqref{eq:complex-affine-lie-brackets}. Set $V=\operatorname{span}\{u,v\}$ and define
  $J:V\to V$ by $Ju=v$, $Jv=-u$. Write $X=at+bj+w$, where $a,b\in\R$ and $w\in V$, and set
  $A=\ad_X$. The brackets give
  \[
    At=-w,\qquad Aj=-Jw,\qquad A|_V=aI+bJ.
  \]
  Since $A$ induces zero on $\g/V$, its spectrum, counted with algebraic multiplicities, is
  \begin{equation}\label{eq:complex-affine-adjoint-spectrum}
    \spec A=\{0,0,a+ib,a-ib\},\qquad \tr A=2a.
  \end{equation}

  Suppose $c\ne0$ and $s\ne0$. By \eqref{eq:trace-and-length-identities}, $s=-mc$. The proof of
  \namedref{Proposition}{prop:even-dimensional-necessity} shows that $-c$ is an eigenvalue of
  $A$. It is real and nonzero, so \eqref{eq:complex-affine-adjoint-spectrum} forces $b=0$ and
  $a=-c$. The trace identity then gives the contradiction
  \[
    -2c=2a=\tr A=-4c-s/m=-3c.
  \]

  Suppose $c\ne0$ and $s=0$. The trace identity gives $2a=-4c$, hence $a=-2c\ne0$. Therefore
  $A|_V=aI+bJ$ is invertible. Since $\im A\subset V$, we have $\rank A=2$ and $\dim\ker A=2$.
  For every $y\in\ker A$, \eqref{eq:trace-and-length-identities} gives
  \[
    -2c\alpha(y)=g(A^*X,y)=g(X,Ay)=0.
  \]
  Thus $\alpha$ vanishes on $\ker A$. Applying \eqref{eq:metric-adjoint-field-equation} to
  $y_1,y_2\in\ker A$ now yields
  \[
    \begin{aligned}
      0&=g(Ay_1,y_2)+g(y_1,Ay_2)\\
      &=-2c\,g(y_1,y_2)-\frac2m\alpha(y_1)\alpha(y_2)
      =-2c\,g(y_1,y_2).
    \end{aligned}
  \]
  Consequently, $\ker A$ is a two-dimensional totally isotropic subspace. This contradicts the
  Lorentzian signature.
\end{proof}

\begin{remark}\label{rem:neutral-metrics-on-lorentzian-exceptions}
  Trans-Lorentzian metrics and corresponding \mMCVF{}s on the two exceptional algebras can be
  written down without further choices. Use $X=u$ and the matrix $G_{(2,2)}$ of
  \eqref{eq:neutral-gram-paired-length-two-chains} in the following ordered bases:
  \[
    \begin{array}{c|c}
      \mathfrak d & (u,-T,z,v)\\
      \aff(\C) & (u,-t,v,-j).
    \end{array}
  \]
  In each case $Ae_2=e_1$ and $Ae_4=e_3$, so $\ad_X$ has Jordan partition $(2,2)$. For this
  partition, \namedref{Proposition}{prop:four-dimensional-jordan-signatures} allows only
  signature $(2,2)$ when $c=0$.
\end{remark}

\Needspace{10\baselineskip}
Part~\textup{(\ref{signature:riemannian-existence})} is the four-dimensional specialization of
\cite[Theorem 4.2]{ZhangChen2023}. We record the condition for each fixed almost-abelian
presentation.

\begin{proof}[Proof of
\namedref{Theorem}{thm:four-dimensional-signature-classification}\textup{(\ref{signature:riemannian-existence})}]
  Suppose a Riemannian metric $g$ admits a nonzero left-invariant \mMCVF{} $X$. By \cite[Theorem
  4.2]{ZhangChen2023},
  \[
    V=[\g,\g]=X^\perp,\qquad \dim V=3,\qquad [V,V]=0.
  \]
  The normalization used there rescales the metric and the field by positive constants, leaving
  $[\g,\g]$ and $X^\perp$ unchanged. Set $A=\ad_X$ and $s=g(X,X)>0$. By
  \namedref{Proposition}{prop:left-invariant-skew-adjoint-reduction}, in a $g$-orthonormal basis
  of $V$ we have
  \[
    c=-s/m<0,\qquad A|_V=-cI_3+S_0,\qquad S_0^T=-S_0.
  \]

  Fix any presentation $\g=\R T\ltimes_D W$ in \namedref{family}{family:almost-abelian}. Since
  $W$ is abelian, $[\g,\g]=DW$. Hence
  \[
    V=DW\subset W,\qquad \dim V=\dim W=3,
  \]
  so $W=V$. Write $X=tT+w$, where $w\in W$. Since $X\notin V$, we have $t\ne0$. As $W$ is
  abelian, $A|_W=tD$, and therefore
  \[
    D\sim_{\R}-\frac ct I_3+\frac1tS_0,\qquad
    \left(\frac1tS_0\right)^T=-\frac1tS_0,\qquad -\frac ct\ne0.
  \]
  This proves necessity for every given presentation.

  Conversely, suppose $D\sim_{\R}aI+S$, with $a\ne0$ and $S^T=-S$. Choose a basis of $\R^3$ in
  which $D=aI+S$. Declare this basis orthonormal and make $T$ a unit vector perpendicular to
  $\R^3$. The same calculation as in \namedref{Example}{ex:odd-dimensional-counterexample}
  applies to every skew-symmetric $S$ and gives
  \[
    X=maT,\qquad c=-ma^2.
  \]
  Thus \eqref{eq:constant-conformal-equation} holds.

  Finally, no algebra in
  \namedref{families}{family:heisenberg-extension}--\ref{family:split-simple-central} has a
  three-dimensional abelian derived algebra. For \namedref{family}{family:heisenberg-extension},
  it is $\im B+\R z$, which is Heisenberg when $B$ is invertible and has dimension at most two
  otherwise. The derived algebras of \namedref{families}{family:split-real-affine}
  and~\ref{family:complex-affine} have dimension two, and that of
  \namedref{family}{family:split-simple-central} is $\slTwoR$. This excludes these families in
  the Riemannian case.
\end{proof}

\namedref{Theorems}{thm:four-dimensional-lie-algebra-classification}
and~\ref{thm:four-dimensional-signature-classification} together prove
\namedref{Theorem}{thm:four-dimensional-existence-and-signatures}.

\bibliographystyle{m_modified_amsplain}
\bibliography{m_modified_references}

@article{AlekseevskyCortesLeistner2023,
  author = {Alekseevsky, D. V. and Cort{\'e}s, V. and Leistner, T.},
  title = {Geometry and holonomy of indecomposable cones},
  journal = {Rev. Mat. Iberoam.},
  volume = {39},
  year = {2023},
  number = {3},
  pages = {1105--1141},
  doi = {10.4171/RMI/1330},
}

@article{AndradaEtAl2005,
  author = {Andrada, A. and Barberis, M. L. and Dotti, I. G. and Ovando, G. P.},
  title = {Product structures on four dimensional solvable {Lie} algebras},
  journal = {Homology Homotopy Appl.},
  volume = {7},
  year = {2005},
  number = {1},
  pages = {9--37},
  archiveprefix = {arXiv},
  eprint = {math/0402234},
  eprintyear = {2004},
}

@misc{AraujoChenLeandro2016,
  author = {Araujo Cintra, A. and Chen, Z. and Leandro Neto, B.},
  title = {Conformal vector fields on {Lie} groups},
  year = {2016},
  archiveprefix = {arXiv},
  eprint = {1608.05943v2},
}

@article{BiggsRemsing2016,
  author = {Biggs, R. and Remsing, C. C.},
  title = {On the classification of real four-dimensional {Lie} groups},
  journal = {J. Lie Theory},
  volume = {26},
  year = {2016},
  number = {4},
  pages = {1001--1035},
  doi = {10.5802/jolt.921},
}

@article{Calvaruso2017,
  author = {Calvaruso, G.},
  title = {Three-dimensional homogeneous generalized {Ricci} solitons},
  journal = {Mediterr. J. Math.},
  volume = {14},
  year = {2017},
  eid = {216},
  pagetotal = {21},
  doi = {10.1007/s00009-017-1019-2},
  archiveprefix = {arXiv},
  eprint = {1503.07767},
  eprintyear = {2015},
}

@article{ChenWolf2018,
  author = {Chen, Z. and Wolf, J. A.},
  title = {Semisimple weakly symmetric pseudo-{Riemannian} manifolds},
  journal = {Abh. Math. Semin. Univ. Hambg.},
  volume = {88},
  year = {2018},
  number = {2},
  pages = {331--369},
  doi = {10.1007/s12188-018-0195-8},
}

@article{CollEtAl2001,
  author = {Coll, B. and Hildebrandt, S. R. and Senovilla, J. M. M.},
  title = {{Kerr--Schild} symmetries},
  journal = {Gen. Relativ. Gravit.},
  volume = {33},
  year = {2001},
  pages = {649--670},
  doi = {10.1023/A:1010265830882},
  archiveprefix = {arXiv},
  eprint = {gr-qc/0006044},
  eprintyear = {2000},
}

@article{DeshmukhEtAl2026ModifiedConformal,
  author = {Deshmukh, S. and Bin Turki, N. and Chen, B.-Y. and Fatima, T.},
  title = {Properties of {$m$}-modified conformal vector fields on a {Riemannian} manifold},
  journal = {AIMS Math.},
  volume = {11},
  year = {2026},
  number = {9},
  pages = {29774--29795},
  doi = {10.3934/math.20261181},
}

@article{DumitrescuZeghib2010,
  author = {Dumitrescu, S. and Zeghib, A.},
  title = {G\'eom\'etries lorentziennes de dimension 3: classification et compl\'etude},
  journal = {Geom. Dedicata},
  volume = {149},
  year = {2010},
  pages = {243--273},
  doi = {10.1007/s10711-010-9480-0},
  archiveprefix = {arXiv},
  eprint = {math/0703846},
  eprintyear = {2007},
}

@article{FreibertSchulteHengesbach2012,
  author = {Freibert, M. and Schulte-Hengesbach, F.},
  title = {Half-flat structures on decomposable {Lie} groups},
  journal = {Transform. Groups},
  volume = {17},
  year = {2012},
  pages = {123--141},
  doi = {10.1007/s00031-011-9168-z},
}

@article{HanounahMehidi2025,
  author = {Hanounah, M. and Mehidi, L.},
  title = {On completeness of foliated structures, and null {Killing} fields},
  journal = {Math. Ann.},
  volume = {392},
  year = {2025},
  pages = {619--658},
  doi = {10.1007/s00208-025-03095-9},
}

@article{JiangYang2025,
  author = {Jiang, J. and Yang, Y.},
  title = {Generalized {Ricci} solitons of four-dimensional non-abelian nilpotent {Lie} groups},
  journal = {J. Nonlinear Math. Phys.},
  volume = {32},
  year = {2025},
  eid = {43},
  pagetotal = {24},
  doi = {10.1007/s44198-025-00296-3},
}

@article{KiosakMatveev2009,
  author = {Kiosak, V. and Matveev, V. S.},
  title = {There are no conformal {Einstein} rescalings of complete pseudo-{Riemannian} {Einstein} metrics},
  journal = {C. R. Math. Acad. Sci. Paris},
  volume = {347},
  year = {2009},
  pages = {1067--1069},
  doi = {10.1016/j.crma.2009.06.017},
}

@book{Knapp2002,
  author = {Knapp, A. W.},
  title = {{Lie} Groups Beyond an Introduction},
  edition = {second},
  series = {Progress in Mathematics},
  volume = {140},
  publisher = {Birkh{\"a}user},
  address = {Boston},
  year = {2002},
  url = {https://www.math.stonybrook.edu/~aknapp/download/Beyond2.pdf},
  urltitle = {Author's digital second edition (2023)},
}

@article{MikesChodorova2010,
  author = {Mike{\v{s}}, J. and Chodorov{\'a}, M.},
  title = {On concircular and torse-forming vector fields on compact manifolds},
  journal = {Acta Math. Acad. Paedagog. Nyh\'azi. (N.S.)},
  volume = {26},
  year = {2010},
  pages = {329--335},
  url = {https://www.emis.de/journals/AMAPN/vol26_2/amapn26_25.pdf},
  urltitle = {Open-access journal text},
}

@article{NurowskiRandall2016,
  author = {Nurowski, P. and Randall, M.},
  title = {Generalized {Ricci} solitons},
  journal = {J. Geom. Anal.},
  volume = {26},
  year = {2016},
  pages = {1280--1345},
  doi = {10.1007/s12220-015-9592-8},
  archiveprefix = {arXiv},
  eprint = {1409.4179},
  eprintyear = {2014},
}

@article{PoddarSharma2025,
  author = {Poddar, R. and Sharma, R.},
  title = {On the triviality of {$m$}-modified conformal vector fields},
  journal = {Indag. Math. (N.S.)},
  volume = {36},
  year = {2025},
  number = {5},
  pages = {1481--1490},
  doi = {10.1016/j.indag.2025.05.009},
  archiveprefix = {arXiv},
  eprint = {2409.07607},
  eprintyear = {2024},
}

@article{ZhangChen2020Conformal,
  author = {Zhang, H. and Chen, Z.},
  title = {Pseudo-{Riemannian} {Lie} groups admitting left-invariant conformal vector fields},
  journal = {C. R. Math. Acad. Sci. Paris},
  volume = {358},
  year = {2020},
  number = {2},
  pages = {143--149},
  doi = {10.5802/crmath.23},
}

@article{ZhangChen2022Lorentzian,
  author = {Zhang, H. and Chen, Z.},
  title = {The classification of {Lorentzian} {Lie} groups with non-{Killing} left-invariant conformal vector fields},
  journal = {Bull. Lond. Math. Soc.},
  volume = {54},
  year = {2022},
  number = {4},
  pages = {1326--1339},
  doi = {10.1112/blms.12631},
}

@article{ZhangChen2023,
  author = {Zhang, H. and Chen, Z.},
  title = {On {$m$}-modified conformal vector fields},
  journal = {J. Geom. Anal.},
  volume = {33},
  year = {2023},
  eid = {258},
  pagetotal = {24},
  doi = {10.1007/s12220-023-01319-5},
}

@article{ZhangChenTan2024,
  author = {Zhang, H. and Chen, Z. and Tan, J.},
  title = {Conformal vector fields on {Lie} groups: The trans-{Lorentzian} signature},
  journal = {J. Algebra},
  volume = {646},
  year = {2024},
  pages = {326--356},
  doi = {10.1016/j.jalgebra.2024.02.010},
}
\end{document}